arXiv:1606.02116v6 [math.OC] 6 Mar 2017

# Local Convergence Properties of Douglas–Rachford and ADMM

Jingwei Liang*, Jalal Fadili* and Gabriel Peyré†

**Abstract.** The Douglas–Rachford (DR) and alternating direction method of multipliers (ADMM) are two proximal splitting algorithms designed to minimize the sum of two proper lower semi-continuous convex functions whose proximity operators are easy to compute. The goal of this work is to understand the local linear convergence behaviour of DR/ADMM when the involved functions are moreover partly smooth. More precisely, when the two functions are partly smooth relative to their respective smooth submanifolds, we show that DR/ADMM (i) identifies these manifolds in finite time; (ii) enters a local linear convergence regime. When both functions are locally polyhedral, we show that the optimal convergence radius is given in terms of the cosine of the Friedrichs angle between the tangent spaces of the identified submanifolds. Under polyhedrality of both functions, we also provide condition sufficient for finite convergence of DR. The obtained results are illustrated by several concrete examples and supported by numerical experiments.

**Key words.** Douglas–Rachford, ADMM, Partial Smoothness, Finite Activity Identification, Local Linear Convergence.

**AMS subject classifications.** 49J52, 65K05, 65K10, 90C25, 90C31.

## 1 Introduction

### 1.1 Non-smooth optimization

In this paper, we consider the structured optimization problem

$$\min_{x\in\mathbb{R}^n} G(x) + J(x), \tag{$\mathcal{P}$}$$

where

- **(A.1)** $G, J \in \Gamma_0(\mathbb{R}^n)$, the class of proper convex and lower semi-continuous (lsc) functions on $\mathbb{R}^n$;
- **(A.2)** $\mathrm{ri}(\mathrm{dom}(G)) \cap \mathrm{ri}(\mathrm{dom}(J)) \neq \emptyset$, where $\mathrm{ri}(C)$ is the relative interior of the nonempty convex set $C$, and $\mathrm{dom}(\cdot)$ denotes the domain of the corresponding function;
- **(A.3)** $\mathrm{Argmin}(G + J) \neq \emptyset$, *i.e.* the set of minimizers is nonempty.

We also assume that these two functions are simple, meaning that their corresponding proximity operators (see (1.2)), $\mathrm{prox}_{\gamma J}$ and $\mathrm{prox}_{\gamma G}$ where $\gamma \in ]0, +\infty[$, are easy to compute, either exactly or up to a very good approximation. Problem ($\mathcal{P}$) covers a large number of problems in areas such as statistical machine learning, inverse problems, signal and image processing to name a few (see Section 9).

### 1.2 Douglas–Rachford splitting

An efficient and provably convergent method to solve problem ($\mathcal{P}$) is the Douglas–Rachford (DR) splitting method. DR was originally proposed in [20] to solve a system of linear equations arising from the discretization of a partial differential equation. The extension of this method suitable to solve optimization and feasibility problems is due to Lions and Mercier [37]. In its exact relaxed form [21, 11, 13], the iteration of DR reads

$$\begin{cases} v_{k+1} = \mathrm{prox}_{\gamma G}(2x_k - z_k), \\ z_{k+1} = (1-\lambda_k)z_k + \lambda_k(z_k + v_{k+1} - x_k), \\ x_{k+1} = \mathrm{prox}_{\gamma J}(z_{k+1}), \end{cases} \tag{1.1}$$

where $\gamma \in ]0, +\infty[$, $\lambda_k \in ]0, 2[$ is the relaxation parameter, and $\mathrm{prox}_{\gamma J} : \mathbb{R}^n \to \mathbb{R}^n$ denotes the proximity operator of $\gamma J$ which is defined as,

$$\mathrm{prox}_{\gamma J}(\cdot) \stackrel{\text{def}}{=} \mathrm{argmin}_{x\in\mathbb{R}^n} \gamma J(x) + \frac{1}{2}\|x - \cdot\|^2. \tag{1.2}$$

*Normandie Univ, ENSICAEN, CNRS, GREYC, E-mail: {Jingwei.Liang, Jalal.Fadili}@ensicaen.fr
†CNRS, DMA, ENS Paris, E-mail: Gabriel.Peyre@ens.fr

In the sequel, we also denote the reflected proximity operator of $\gamma J$ as

$$\mathrm{rprox}_{\gamma J} \stackrel{\text{def}}{=} 2\mathrm{prox}_{\gamma J} - \mathrm{Id}.$$

The DR scheme (1.1) can be rewritten into a fixed-point iteration with respect to $\{z_k\}_{k\in\mathbb{N}}$, *i.e.*

$$z_{k+1} = \mathscr{F}_{\gamma,\lambda_k}(z_k), \tag{1.3}$$

where

$$\mathscr{F}_{\gamma,\lambda_k} \stackrel{\text{def}}{=} (1-\lambda_k)\mathrm{Id} + \lambda_k \mathscr{F}_\gamma \quad \text{and} \quad \mathscr{F}_\gamma \stackrel{\text{def}}{=} \frac{\mathrm{rprox}_{\gamma G} \circ \mathrm{rprox}_{\gamma J} + \mathrm{Id}}{2}. \tag{1.4}$$

Under assumptions (**A.1**)-(**A.3**), and if $\sum_{k\in\mathbb{N}} \lambda_k(2-\lambda_k) = +\infty$, it is known that $z_k$ converges to some fixed point $z^\star \in \mathrm{Fix}(\mathscr{F}_\gamma) \neq \emptyset$, and that the shadow point $x_k$ and $v_k$ both converge to $x^\star \stackrel{\text{def}}{=} \mathrm{prox}_{\gamma J}(z^\star) \in \mathrm{Argmin}(G + J)$; see *e.g.* [4, Corollary 27.7].

In this paper, we consider a non-stationary version of (1.1) (*i.e.* the choices of $\gamma$ vary along the iterations), which is described below in Algorithm 1.

---

**Algorithm 1:** Nonstationary Douglas–Rachford splitting

---

**Initial**: $k = 0$, $z_0 \in \mathbb{R}^n$, $x_0 = \mathrm{prox}_{\gamma_0 J}(z_0)$;

**repeat**

Let $\gamma_k \in ]0, +\infty[$, $\lambda_k \in ]0, 2[$:

$$\begin{aligned} v_{k+1} &= \mathrm{prox}_{\gamma_k G}(2x_k - z_k), \\ z_{k+1} &= (1-\lambda_k)z_k + \lambda_k(z_k + v_{k+1} - x_k), \\ x_{k+1} &= \mathrm{prox}_{\gamma_{k+1} J}(z_{k+1}), \end{aligned} \tag{1.5}$$

$k = k + 1$;

**until** *convergence*;

---

**Remark 1.1.**

(i) By definition, the DR method is not symmetric with respect to the order of the functions $J$ and $G$, see [6] for a systematic study of the two possible versions in the exact, stationary and unrelaxed case. Nevertheless, all of our statements throughout hold true, with minor adaptations, when the order of $J$ and $G$ is reversed in (1.5). Note also that the standard DR only accounts for the sum of two functions. Extension to more than two functions is straightforward through a product space trick, see Section 8 for details.

(ii) This paper consists of two main parts, the first one dealing with global convergence guarantees of (1.5) (Section 3), and a second one on its the local convergence properties when the involved functions are also partly smooth (Section 5). It is for the sake of the latter that we mainly focus on the finite dimensional setting $\mathbb{R}^n$. It is worth pointing out, however, that the global convergence result (Theorem 3.1) also holds for real Hilbert space case where weak convergence can be obtained.

(iii) For global convergence, one can also consider an inexact version of (1.5) by incorporating additive errors in the computation of $x_k$ and $v_k$, though we do not elaborate more on this for the sake of local convergence analysis.

### 1.3 Contributions

Global sublinear convergence rate estimates of DR and the alternating direction method of multipliers (ADMM) iterations have been recently established in the literature, see *e.g.* [35, 17] and the references therein. Such rate becomes linear under further assumptions, typically smoothness and strong convexity, see e.g. [18, 23] and references therein. DR method however shows local linear convergence in various situations where strong convexity is absent. Studying the local linear convergence of DR or ADMM without strong convexity has received attention in recent years, see the detailed discussion in Section 1.4.

However, most of the existing work either focuses on some special cases where a specific structure of the problem at hand can be exploited, or imposes certain regularity conditions which are barely verified in practical situations. Therefore, it is important to present a unified analysis framework, and possibly with stronger claims. This is one of the main motivations of this work. More precisely, our main contributions are the following.

**Globally convergent non-stationary DR** By casting the non-stationarity as an additional error, in Section 3 we establish a global convergence result for the non-stationary DR iteration (1.5). This turns exploits our previous result on the convergence of the general inexact and non-stationary Krasnosel'skiĭ–Mann fixed-point iteration studied in [35].

**Finite time activity identification** Let $x^\star \in \mathrm{Argmin}(G+J) \neq \emptyset$. Then assuming, in addition to (**A.1**)-(**A.3**), that both $G$ and $J$ are partly smooth at $x^\star$ relative to $C^2$-smooth submanifolds $\mathcal{M}^G_{x^\star}$ and $\mathcal{M}^J_{x^\star}$, respectively, (see Definition 4.1), we show in Section 4.1 that under a non-degeneracy condition, the non-stationary DR sequences $\{v_k\}_{k\in\mathbb{N}}$ and $\{x_k\}_{k\in\mathbb{N}}$ respectively identify in finite time the submanifolds $\mathcal{M}^G_{x^\star}$ and $\mathcal{M}^J_{x^\star}$. In plain words, this means that after a finite number of iterations, say $K$, we have $x_k \in \mathcal{M}^J_{x^\star}$ and $v_k \in \mathcal{M}^G_{x^\star}$ for all $k \geq K$.

**Local linear convergence** Exploiting the finite identification property, we then show that the non-stationary DR iterates converge locally linearly in Section 5. We characterize the convergence rate precisely based on the properties of the identified partial smoothness submanifolds. Moreover, when both $G$ and $J$ are locally polyhedral around $x^\star$ and the stationary DR scheme is considered with constant $\gamma_k$ and $\lambda_k$, we show that the optimal convergence rate is given in terms of the *cosine* of the *Friedrichs angle* between the tangent spaces of the two submanifolds. We also generalize these claims to the minimization of the sum of more than two functions.

**Finite convergence** Building upon our local convergence analysis, we also characterize situations where finite convergence occurs in Section 6. More precisely, when the stationary and unrelaxed DR iteration is used and $G$ and $J$ are locally polyhedral nearby $x^\star$, and if either $G$ or $J$ is differentiable at $x^\star$, we obtain finite convergence.

We also touch one some practical acceleration schemes, since once the active submanifolds are identified, the globally convex but non-smooth problem becomes locally $C^2$-smooth, though possibly non-convex. As a consequence, it opens the door to high-order optimization methods, such as Newton-like or nonlinear conjugate gradient.

**ADMM** It have been well-known the relation between alternating direction method of multipliers and DR, that ADMM is applying DR to the Fenchel dual of ($\mathcal{P}$) (see [22, 21]). By virtue of such relation, we can adapt our obtained result of DR with proper modifications. Therefore, in Section 7, we present in parallel the finite activity identification and local linear convergence result for the ADMM.

### 1.4 Relation to prior work

There are problem instances in the literature where the (stationary) DR and ADMM algorithms are proved to converge linearly either globally or locally. For instance, in [37, Proposition 4], it is assumed that the "internal" function is strongly convex with a Lipschitz continuous gradient. This local linear convergence result is further investigated in [17, 23] under smoothness and strong convexity assumptions. The special case of Basis Pursuit (BP), *i.e.* $\ell_1$ minimization with an affine constraint, is considered in [19] and an eventual local linear convergence is shown in the absence of strong convexity. In [7], the author analyses the local convergence behaviour of ADMM for quadratic or linear programming, and shows local linear convergence if the optimal solution is unique and strict complementarity holds. For the case of two subspaces (though in general real Hilbert space), linear convergence of DR with the optimal rate being the cosine of the Friedrichs angle between the subspaces is proved in [2]. It turns out that [19, 7, 2] are special cases of our framework, and our results generalize theirs to a larger class of problems. The proposed work is also a more general extension of our previous results in [36] which tackled only the case of locally polyhedral functions.

For the non-convex case, [8] considers DR for a feasibility problem of a sphere intersecting a line or more generally a proper affine subset. Such feasibility problems with an affine subspace and a super-regular set (in the sense of [32]) with strongly regular intersection is considered in [27], and is generalized later to two $(\varepsilon, \delta)$-regular sets with linearly regular intersection [26], see also [40] for an even more general setting. However, even in the convex polyhedral case, the rate provided in [40].

Our finite convergence result complements and extends that of [5] who established finite convergence of (unrelaxed stationary) DR in the presence of Slater's condition, for solving convex feasibility problems where one set is an affine subspace and the other is a polyhedron.

**Paper organization** The rest of the paper is organized as follows. Some ingredients from monotone operator theory and angles between subspaces are collected in Section 2. Global convergence of the non-stationary DR (1.5) is established in Section 3. Section 4 discusses the notion of partial smoothness and some essential properties. We then turn to the main contributions of this paper, namely finite time activity identification (Section 4.1), local linear convergence (Section 5) and finite termination (Section 6) of DR under partial smoothness. We develop in parallel the result for ADMM under partial smoothness in Section 7. Section 8 extends the results to the sum of more than two functions. In Section 9, we report various numerical experiments to support our theoretical findings.

## 2 Preliminaries

Throughout the paper, $\mathbb{N}$ is the set of nonnegative integers, $\mathbb{R}^n$ is a $n$-dimensional real Euclidean space equipped with scalar product $\langle \cdot, \cdot \rangle$ and norm $\|\cdot\|$. $\mathrm{Id}$ denotes the identity operator on $\mathbb{R}^n$. For a vector $x \in \mathbb{R}^n$ and some subset of indices $b \subset \{1, \ldots, n\}$, $x_b$ is the restriction of $x$ to the entries indexed in $b$. For $p \in [1, +\infty]$, $\|x\|_p \stackrel{\text{def}}{=} (\sum_{i=1}^p |x_i|)^{1/p}$ is the $\ell_p$ in $\mathbb{R}^n$ with the usual adaptation for $p = +\infty$. $\ell_+^1$ denotes the set of summable sequences in $[0, +\infty[$. For a matrix $M \in \mathbb{R}^{n \times n}$, we denote $\|M\|$ its operator norm and $\rho(M)$ its spectral radius.

The subdifferential of a function $J \in \Gamma_0(\mathbb{R}^n)$ is the set-valued operator,

$$\partial J : \mathbb{R}^n \rightrightarrows \mathbb{R}^n, \, x \mapsto \big\{ g \in \mathbb{R}^n | J(y) \geq J(x) + \langle g, \, y - x \rangle, \forall y \in \mathbb{R}^n \big\}.$$

Recall the proximity operator defined in (1.2), it can be verified that it is also the *resolvent* of the subdifferential, *i.e.* $\mathrm{prox}_{\gamma J} = (\mathrm{Id} + \gamma\partial J)^{-1}$.

For a nonempty convex set $C \subset \mathbb{R}^n$, denote $\mathrm{cone}(C)$ its conical hull, $\mathrm{aff}(C)$ its affine hull, and $\mathrm{par}(C) = \mathbb{R}(C - C)$ the subspace parallel to $C$, *i.e.* a translate of $\mathrm{aff}(C)$ to the origin. $\mathrm{P}_C$ is the orthogonal projection operator onto $C$ and $N_C(x)$ its normal cone at $x$.

## 2.1 Operators and matrices

**Definition 2.1 (Monotone operator).** A set-valued operator $\mathscr{A} : \mathbb{R}^n \rightrightarrows \mathbb{R}^n$ is monotone if, given any $x, z \in \mathbb{R}^n$, there holds

$$\langle x - z,\, u - v\rangle \geq 0,\ \forall (x, u) \in \mathrm{gph}(\mathscr{A}) \ \text{ and } \ (z, v) \in \mathrm{gph}(\mathscr{A}),$$

where $\mathrm{gph}(\mathscr{A}) \stackrel{\mathrm{def}}{=} \big\{(x, u) \in \mathbb{R}^n \times \mathbb{R}^n : \ u \in \mathscr{A}(x)\big\}$. It is moreover maximal monotone if its graph is not strictly contained in the graph of any other monotone operators.

The best-known example of maximal monotone operator is the subdifferential mapping of functions in $\Gamma_0(\mathbb{R}^n)$.

**Definition 2.2 (Non-expansive operator).** An operator $\mathscr{F} : \mathbb{R}^n \to \mathbb{R}^n$ is nonexpansive if

$$\forall x, y \in \mathbb{R}^n,\ \ \|\mathscr{F}(x) - \mathscr{F}(y)\| \leq \|x - y\|.$$

For any $\alpha \in ]0, 1[$, $\mathscr{F}$ is called $\alpha$-averaged if there exists a nonexpansive operator $\mathscr{R}$ such that $\mathscr{F} = \alpha\mathscr{R} + (1 - \alpha)\mathrm{Id}$.

The class of $\alpha$-averaged operators is closed under relaxation, convex combination and composition [4, 15]. In particular when $\alpha = \frac{1}{2}$, $\mathscr{F}$ is called *firmly nonexpansive*. Several properties of firmly nonexpansive operators are collected in the following lemma.

**Lemma 2.3.** *Let $\mathscr{F} : \mathbb{R}^n \to \mathbb{R}^n$. Then the following statements are equivalent:*

(i) *$\mathscr{F}$ is firmly nonexpansive;*

(ii) *$\mathrm{Id} - \mathscr{F}$ is firmly nonexpansive;*

(iii) *$2\mathscr{F} - \mathrm{Id}$ is nonexpansive;*

(iv) *Given any $\lambda \in ]0, 2]$, $(1 - \lambda)\mathrm{Id} + \lambda\mathscr{F}$ is $\frac{\lambda}{2}$-averaged;*

(v) *$\mathscr{F}$ is the resolvent of a maximal monotone operator $A : \mathbb{R}^n \rightrightarrows \mathbb{R}^n$.*

**Proof.** (i)⇔(ii)⇔(iii) follow from [4, Proposition 4.2, Corollary 4.29], (i)⇔(iv) is [4, Corollary 4.29], and (i)⇔(v) is [4, Corollary 23.8]. □

Recall the fixed-point operator $\mathscr{F}_\gamma$ from (1.4).

**Lemma 2.4.** *$\mathscr{F}_\gamma$ is firmly nonexpansive.*

**Proof.** The reflected resolvents $\mathrm{rprox}_{\gamma J}$ and $\mathrm{rprox}_{\gamma G}$ are nonexpansive [4, Corollary 23.10(ii)], and so is their composition. The claim then follows from Lemma 2.3(i)⇔(ii). □

**Definition 2.5 (Convergent matrices).** A matrix $M \in \mathbb{R}^{n\times n}$ is convergent to $M^\infty \in \mathbb{R}^{n\times n}$ if, and only if,

$$\lim_{k\to\infty} \|M^k - M^\infty\| = 0.$$

The following identity is known as the spectral radius formula

$$\rho(M) = \lim_{k\to+\infty} \|M^k\|^{1/k}. \tag{2.1}$$

## 2.2 Angles between subspaces

In this part we introduce the principal angles and the Friedrichs angle between two subspaces $T_1$ and $T_2$. Without loss of generality, let $1 \leq p \stackrel{\text{def}}{=} \dim(T_1) \leq q \stackrel{\text{def}}{=} \dim(T_2) \leq n-1$.

**Definition 2.6 (Principal angles).** The principal angles $\theta_k \in [0, \frac{\pi}{2}]$, $k = 1, \dots, p$ between subspaces $T_1$ and $T_2$ are defined by, with $u_0 = v_0 \stackrel{\text{def}}{=} 0$

$$\cos(\theta_k) \stackrel{\text{def}}{=} \langle u_k,\, v_k\rangle = \max\langle u,\, v\rangle \ \text{ s.t. } \ u \in T_1, v \in T_2, \|u\| = 1, \|v\| = 1,$$
$$\langle u,\, u_i\rangle = \langle v,\, v_i\rangle = 0,\ i = 0, \cdots, k-1.$$

The principal angles $\theta_k$ are unique with $0 \leq \theta_1 \leq \theta_2 \leq \cdots \leq \theta_p \leq \pi/2$.

**Definition 2.7 (Friedrichs angle).** The Friedrichs angle $\theta_F \in ]0, \frac{\pi}{2}]$ between $T_1$ and $T_2$ is

$$\cos\big(\theta_F(T_1, T_2)\big) \stackrel{\text{def}}{=} \max\langle u,\, v\rangle \ \text{ s.t. } \ u \in T_1 \cap (T_1 \cap T_2)^{\perp}, \|u\| = 1,\ v \in T_2 \cap (T_1 \cap T_2)^{\perp}, \|v\| = 1.$$

The following lemma shows the relation between the Friedrichs and principal angles whose proof can be found in [3, Proposition 3.3].

**Lemma 2.8 (Principal angles and Friedrichs angle).** *The Friedrichs angle is exactly $\theta_{d+1}$ where $d \stackrel{\text{def}}{=} \dim(T_1 \cap T_2)$. Moreover, $\theta_F(T_1, T_2) > 0$.*

**Remark 2.9.** One approach to obtain the principal angles is through the singular value decomposition (SVD). For instance, let $X \in \mathbb{R}^{n\times p}$ and $Y \in \mathbb{R}^{n\times q}$ form orthonormal bases for the subspaces $T_1$ and $T_2$ respectively. Let $U\Sigma V^T$ be the SVD of $X^TY \in \mathbb{R}^{p\times q}$, then $\cos(\theta_k) = \sigma_k$, $k = 1, 2, \dots, p$ and $\sigma_k$ corresponds to the $k$'th largest singular value in $\Sigma$.

# 3 Global convergence

Recall the operators defined in (1.4). The nonstationay DR iteration (1.5) can also be written

$$z_{k+1} = \mathscr{F}_{\gamma_k,\lambda_k}(z_k) = \mathscr{F}_{\gamma,\lambda_k}(z_k) + (\mathscr{F}_{\gamma_k,\lambda_k} - \mathscr{F}_{\gamma,\lambda_k})(z_k). \tag{3.1}$$

In plain words, the nonstationary iteration (1.5) can be seen as a perturbed version of the stationary one (1.1).

**Theorem 3.1 (Global convergence).** *Consider the nonstationary DR iteration (1.5). Suppose that the following conditions are fulfilled*

**(H.1)** *Assumptions (**A.1**)-(**A.3**) hold;*

**(H.2)** $\lambda_k \in [0, 2]$ *such that* $\sum_{k\in\mathbb{N}} \lambda_k(2-\lambda_k) = +\infty$*;*

**(H.3)** $(\gamma_k, \gamma) \in ]0, +\infty[^2$ *such that*
$$\{\lambda_k|\gamma_k - \gamma|\}_{k\in\mathbb{N}} \in \ell^1_+.$$

*Then the sequence $\{z_k\}_{k\in\mathbb{N}}$ converges to a point $z^\star \in \mathrm{Fix}(\mathscr{F}_\gamma)$ and $x^\star = \mathrm{prox}_{\gamma J}(z^\star) \in \mathrm{Argmin}(G{+}J)$. Moreover, the shadow sequence $\{x_k\}_{k\in\mathbb{N}}$ and $\{v_k\}_{k\in\mathbb{N}}$ both converge to $x^\star$ if $\gamma_k$ is convergent.*

See Section A for the proof.

**Remark 3.2.**

(i) The conclusions of Theorem 3.1 remain true if $x_k$ and $v_k$ are computed inexactly with additive errors $\varepsilon_{1,k}$ and $\varepsilon_{2,k}$, provided that $\{\lambda_k\|\varepsilon_{1,k}\|\}_{k\in\mathbb{N}} \in \ell^1_+$ and $\{\lambda_k\|\varepsilon_{2,k}\|\}_{k\in\mathbb{N}} \in \ell^1_+$.

(ii) The summability assumption (**H.3**) is weaker than imposing it without $\lambda_k$. Indeed, following the discussion in [12, Remark 5.7], take $q \in ]0,1]$, and let

$$\lambda_k = 1 - \sqrt{1 - 1/k} \text{ and } |\gamma_k - \gamma| = \frac{1 + \sqrt{1 - 1/k}}{k^q},$$

then it can be verified that

$$|\gamma_k - \gamma| \notin \ell^1_+, \; \lambda_k |\gamma_k - \gamma| = \frac{1}{k^{1+q}} \in \ell^1_+ \text{ and } \lambda_k(2 - \lambda_k) = \frac{1}{k} \notin \ell^1_+.$$

(iii) The assumptions made on the sequence $\{\gamma_k\}_{k\in\mathbb{N}}$ imply that $\gamma_k \to \gamma$ (see Lemma A.1). In fact, if $\inf_{k\in\mathbb{N}} \lambda_k > 0$, we have $\{|\gamma_k - \gamma|\}_{k\in\mathbb{N}} \in \ell_1^+$, entailing $\gamma_k \to \gamma$, and thus the convergence assumption on $\gamma_k$ is superfluous.

# 4 Partial smoothness

The concept of partial smoothness was formalized in [31]. This notion, as well as that of identifiable surfaces [47], captures the essential features of the geometry of non-smoothness which are along the so-called active/identifiable submanifold. For convex functions, a closely related idea is developed in [30]. Loosely speaking, a partly smooth function behaves smoothly as we move along the identifiable submanifold, and sharply if we move transversal to the manifold. In fact, the behaviour of the function and of its minimizers depend essentially on its restriction to this manifold, hence offering a powerful framework for algorithmic and sensitivity analysis theory.

Let $\mathcal{M}$ be a $C^2$-smooth embedded submanifold of $\mathbb{R}^n$ around a point $x$. To lighten terminology, henceforth we shall state $C^2$-manifold instead of $C^2$-smooth embedded submanifold of $\mathbb{R}^n$. The natural embedding of a submanifold $\mathcal{M}$ into $\mathbb{R}^n$ permits to define a Riemannian structure on $\mathcal{M}$, and we simply say $\mathcal{M}$ is a Riemannian manifold. $\mathcal{T}_{\mathcal{M}}(x)$ denotes the tangent space to $\mathcal{M}$ at any point near $x$ in $\mathcal{M}$. More material on manifolds is given in Section C.1.

We are now in position to formally define the class of partly smooth functions in $\Gamma_0(\mathbb{R}^n)$.

**Definition 4.1 (Partly smooth function).** Let $F \in \Gamma_0(\mathbb{R}^n)$, and $x \in \mathbb{R}^n$ such that $\partial F(x) \neq \emptyset$. $F$ is then said to be *partly smooth* at $x$ relative to a set $\mathcal{M}$ containing $x$ if

(i) **Smoothness**: $\mathcal{M}$ is a $C^2$-manifold around $x$, $F$ restricted to $\mathcal{M}$ is $C^2$ around $x$;

(ii) **Sharpness**: The tangent space $\mathcal{T}_{\mathcal{M}}(x)$ coincides with $T_x = \mathrm{par}\big(\partial F(x)\big)^{\perp}$;

(iii) **Continuity**: The set-valued mapping $\partial F$ is continuous at $x$ relative to $\mathcal{M}$.

The class of partly smooth functions at $x$ relative to $\mathcal{M}$ is denoted as $\mathrm{PSF}_x(\mathcal{M})$.

In fact, local polyhedrality also implies that the subdifferential is locally constant around $x$ along $x + T_x$. Capitalizing on the results of [31], it can be shown that under mild transversality assumptions, the set of partly smooth functions is closed under addition and pre-composition by a linear operator. Moreover, absolutely permutation-invariant convex and partly smooth functions of the singular values of a real matrix, *i.e.* spectral functions, are convex and partly smooth spectral functions of the matrix [16]. Some examples of partly smooth functions will be discussed in Section 9.

The next lemma gives expressions of the Riemannian gradient and Hessian (see Section C.1 for definitions) of a partly smooth function.

**Lemma 4.2.** *If $F \in \mathrm{PSF}_x(\mathcal{M})$, then for any $x' \in \mathcal{M}$ near $x$*

$$\nabla_{\mathcal{M}} F(x') = \mathrm{P}_{T_{x'}}(\partial F(x')).$$

*In turn, for all $h \in T_{x'}$*

$$\nabla^2_{\mathcal{M}} F(x')h = \mathrm{P}_{T_{x'}} \nabla^2 \widetilde{F}(x')h + \mathfrak{W}_{x'}\big(h, \mathrm{P}_{T_{x'}^{\perp}} \nabla \widetilde{F}(x')\big),$$

*where $\widetilde{F}$ is any smooth extension (representative) of $F$ on $\mathcal{M}$, and $\mathfrak{W}_x(\cdot,\cdot) : T_x \times T_x^{\perp} \to T_x$ is the Weingarten map of $\mathcal{M}$ at $x$.*

**Proof.** See [34, Fact 3.3]. □

## 4.1 Finite activity identification

With the above global convergence result at hand, we are now ready to state the finite time activity identification property of the non-stationary DR method.

Let $z^\star \in \mathrm{Fix}(\mathscr{F}_{\gamma,\lambda})$ and $x^\star = \mathrm{prox}_{\gamma J}(z^\star) \in \mathrm{Argmin}(G+J)$, then at convergence of the DR iteration (1.5), we have the following inclusion holds,

$$x^\star - z^\star \in \gamma\partial G(x^\star) \ \text{ and } \ z^\star - x^\star \in \gamma\partial J(x^\star).$$

The condition needed for identification result is built upon this inclusion.

**Theorem 4.3 (Finite activity identification).** *For the DR iteration (1.5), suppose that (**H.1**)-(**H.3**) hold and $\gamma_k$ is convergent, entailing that $(z_k, x_k, v_k) \to (z^\star, x^\star, x^\star)$, where $z^\star \in \mathrm{Fix}(\mathscr{F}_{\gamma,\lambda})$ and $x^\star \in \mathrm{Argmin}(G+J)$. Assume also that $\inf_{k\in\mathbb{N}} \gamma_k \geq \underline{\gamma} > 0$. If $G \in \mathrm{PSF}_{x^\star}(\mathcal{M}^G_{x^\star})$ and $J \in \mathrm{PSF}_{x^\star}(\mathcal{M}^J_{x^\star})$, and the non-degeneracy condition*

$$x^\star - z^\star \in \gamma\mathrm{ri}\big(\partial G(x^\star)\big) \ \text{ and } \ z^\star - x^\star \in \gamma\mathrm{ri}\big(\partial J(x^\star)\big) \tag{ND}$$

*holds. Then*

(i) *There exists $K \in \mathbb{N}$ large enough such that for all $k \geq K$, $(v_k, x_k) \in \mathcal{M}^G_{x^\star} \times \mathcal{M}^J_{x^\star}$.*

(ii) *Moreover,*
 (a) *If $\mathcal{M}^J_{x^\star} = x^\star + T^J_{x^\star}$, then $\forall k \geq K$, $T^J_{x_k} = T^J_{x^\star}$.*
 (b) *If $\mathcal{M}^G_{x^\star} = x^\star + T^G_{x^\star}$, then $\forall k \geq K$, $T^G_{v_k} = T^G_{x^\star}$.*
 (c) *$J$ is locally polyhedral around $x^\star$, then $\forall k \geq K$, $x_k \in \mathcal{M}^J_{x^\star} = x^\star + T^J_{x^\star}$, $T^J_{x_k} = T^J_{x^\star}$, $\nabla_{\mathcal{M}^J_{x^\star}} J(x_k) = \nabla_{\mathcal{M}^J_{x^\star}} J(x^\star)$, and $\nabla^2_{\mathcal{M}^J_{x^\star}} J(x_k) = 0$.*
 (d) *$G$ is locally polyhedral around $x^\star$, then $\forall k \geq K$, $v_k \in \mathcal{M}^G_{x^\star} = x^\star + T^G_{x^\star}$, $T^G_{v_k} = T^G_{x^\star}$, $\nabla_{\mathcal{M}^G_{x^\star}} G(v_k) = \nabla_{\mathcal{M}^G_{x^\star}} G(x^\star)$, and $\nabla^2_{\mathcal{M}^G_{x^\star}} G(v_k) = 0$.*

See Section B.1 for the proof.

**Remark 4.4.**

(i) The theorem remains true if the condition on $\gamma_k$ is replaced with $\gamma_k \geq \underline{\gamma} > 0$ and $\lambda_k \geq \underline{\lambda} > 0$, (use (**H.3**) in the proof).

(ii) The non-degeneracy condition (ND) is a geometric generalization of strict complementarity in nonlinear programming. Building on the arguments of [25], it is almost a necessary condition for the finite identification of $\mathcal{M}_{x^\star}$. Relaxing it in general is a challenging problem.

(iii) In general, we have no identification guarantees for $x_k$ and $v_k$ if the proximity operators are computed with errors, even if the latter are summable, in which case one can still prove global convergence (see Remark 3.2). The deep reason behind this is that in the exact case, under condition (ND), the proximal mapping of a partly smooth function and that of its restriction to the corresponding active manifold locally agree nearby $x^\star$. This property can be easily violated if approximate proximal mappings are involved.

(iv) When the minimizer is unique, using the fixed-point set characterization of DR, see e.g. [13, Lemma 2.6], it can be shown that condition (ND) is also equivalent to $z^\star \in \mathrm{ri}(\mathrm{Fix}(\mathscr{F}_\gamma))$.

**A bound on the finite identification iteration** In Theorem 4.3, we have not provided an estimate of the number of iterations beyond which finite identification occurs. In fact, there is a situation where the answer is trivial, *i.e.* $J$ (resp. $G$) is the indicator function of a subspace. However, answering such a question in general remains challenging. In the following, we shall give a bound in some important cases.

We start with the following general statement and then show that it holds true typically for indicators of polyhedral sets. Denote $\tau_k \stackrel{\text{def}}{=} \lambda_k(2-\lambda_k)$.

**Proposition 4.5.** *Suppose that the assumptions of Theorem 4.3 hold, that* $\inf_{k\in\mathbb{N}} \tau_k \geq \underline{\tau} > 0$*, and moreover, the iterates are such that* $\partial J(x_k) \subset \mathrm{rbd}(\partial J(x^\star))$ *whenever* $x_k \notin \mathcal{M}^J_{x^\star}$ *and* $\partial G(v_k) \subset \mathrm{rbd}(\partial G(x^\star))$ *whenever* $v_k \notin \mathcal{M}^G_{x^\star}$*. Then,* $\mathcal{M}^J_{x^\star}$ *and* $\mathcal{M}^G_{x^\star}$ *will be identified for some* $k$ *obeying*

$$k \geq \frac{\|z_0 - z^\star\|^2 + O(\sum_{k\in\mathbb{N}} \lambda_k|\gamma_k - \gamma|)}{\underline{\gamma}^2 \underline{\tau} \mathrm{dist}(0, \mathrm{rbd}(\partial J(x^\star) + \partial G(x^\star)))^2}. \tag{4.1}$$

See Section B.2 for the proof.

Observe that the assumption on $\tau_k$ automatically implies (**H.2**). As one intuitively expects, this upper-bound increases as (ND) becomes more demanding.

**Example 4.6 (Indicators of polyhedral sets).** We will discuss the case of $J$, and the same reasoning applies to $G$. Consider $J$ as the indicator function of a polyhedral set $C^J$, *i.e.*

$$J(x) = \iota_{C^J}(x), \quad \text{where} \quad C^J = \big\{x \in \mathbb{R}^n : \langle c_i^J, x\rangle \leq d_i^J, i = 1, \ldots, m\big\}.$$

Define $I_x^J \stackrel{\text{def}}{=} \big\{i : \langle c_i, x\rangle = d_i\big\}$ the active set at $x$. The normal cone to $C^J$ at $x \in C^J$ is polyhedral and given by [44, Theorem 6.46]

$$\partial J(x) = N_{C^J}(x) = \mathrm{cone}((c_i^J)_{i\in I_x^J}).$$

It is immediate then to show that $J$ is partly smooth at $x \in C^J$ relative to the affine subspace $\mathcal{M}^J_x = x + T_x^J$, where, $T_x^J = \mathrm{span}((c_i^J)_{i\in I_x^J})^\perp$. Let $\mathcal{F}^J_x$ be the face of $C^J$ containing $x$. From [43, Theorem 18.8], one can deduce that

$$\mathcal{F}^J_x = \mathcal{M}^J_x \cap C^J. \tag{4.2}$$

We then have

$$\mathcal{M}^J_{x^\star} \subsetneq \mathcal{M}^J_x \underset{(4.2)}{\Longleftrightarrow} \mathcal{F}^J_{x^\star} \subsetneq \mathcal{F}^J_x \tag{4.3}$$

$$\underset{[28,\ \text{Proposition 3.4}]}{\Longrightarrow} N_{C^J}(x) \text{ is a face of (other than) } N_{C^J}(x^\star) \underset{[43,\ \text{Corollary 18.1.3}]}{\Longrightarrow} \partial J(x) \subset \mathrm{rbd}(\partial J(x^\star)). \tag{4.4}$$

Suppose that $\mathcal{M}^J_{x^\star}$ has not been identified yet. Therefore, since $x_k = \mathrm{P}_{C^J}(z_k) = \mathrm{P}_{\mathcal{F}^J_{x_k} \setminus \mathcal{F}^J_{x^\star}}(z_k)$, and thanks to (4.3), this is equivalent to

$$\text{either} \quad \mathcal{F}^J_{x^\star} \subsetneq \mathcal{F}^J_{x_k} \quad \text{or} \quad \mathcal{F}^J_{x_k} \cap \mathcal{F}^J_{x^\star} = \emptyset.$$

It then follows from (4.4) and Proposition 4.5 that the number of iterations where $\mathcal{F}^J_{x^\star} \subsetneq \mathcal{F}^J_{x_k}$ and $\mathcal{F}^G_{x^\star} \subsetneq \mathcal{F}^G_{x_k}$ cannot exceed the bound in (4.1), and thus identification will happen indeed for some large enough $k$ obeying (4.1).

# 5 Local linear convergence

Building upon the identification result from the previous section, we now turn to the local behaviour of the DR iteration (3.1) under partial smoothness. The key feature is that, once the active manifolds are identified, the DR iteration locally linearizes (possibly up to first-order). It is then sufficient to control the spectral properties of the matrix appearing in the linearized iteration to exhibit the local linear convergence rate.

## 5.1 Locally linearized iteration

Let $z^\star \in \mathrm{Fix}(\mathscr{F}_{\gamma,\lambda})$ and $x^\star = \mathrm{prox}_{\gamma J}(z^\star) \in \mathrm{Argmin}(G+J)$. Define the following two functions

$$\overline{G}(x) \stackrel{\text{def}}{=} \gamma G(x) - \langle x,\, x^\star - z^\star\rangle,\;\; \overline{J}(x) \stackrel{\text{def}}{=} \gamma J(x) - \langle x,\, z^\star - x^\star\rangle. \tag{5.1}$$

We start with the following key lemma.

**Lemma 5.1.** *Suppose that $G \in \mathrm{PSF}_{x^\star}(\mathcal{M}^G_{x^\star})$ and $J \in \mathrm{PSF}_{x^\star}(\mathcal{M}^J_{x^\star})$. Define the two matrices*

$$H_{\overline{G}} \stackrel{\text{def}}{=} \mathrm{P}_{T^G_{x^\star}} \nabla^2_{\mathcal{M}^G_{x^\star}} \overline{G}(x^\star)\mathrm{P}_{T^G_{x^\star}} \;\; \textit{and} \;\; H_{\overline{J}} \stackrel{\text{def}}{=} \mathrm{P}_{T^J_{x^\star}} \nabla^2_{\mathcal{M}^J_{x^\star}} \overline{J}(x^\star)\mathrm{P}_{T^J_{x^\star}}. \tag{5.2}$$

*Then $H_{\overline{G}}$ and $H_{\overline{J}}$ are symmetric positive semi-definite under either of the following circumstances:*
(i) (ND) *holds.*
(ii) $\mathcal{M}^G_{x^\star}$ *and* $\mathcal{M}^J_{x^\star}$ *are affine subspaces.*
*In turn, the matrices*

$$W_{\overline{G}} \stackrel{\text{def}}{=} (\mathrm{Id} + H_{\overline{G}})^{-1} \;\; \textit{and} \;\; W_{\overline{J}} \stackrel{\text{def}}{=} (\mathrm{Id} + H_{\overline{J}})^{-1}, \tag{5.3}$$

*are both firmly non-expansive.*

**Proof.** Here we prove the case for $J$ since the same arguments apply to $G$ just as well. Claims (i) and (ii) follow from [34, Lemma 4.3] since $J \in \mathrm{PSF}_{x^\star}(\mathcal{M}^J_{x^\star})$. Consequently, $W_{\overline{J}}$ is symmetric positive definite with eigenvalues in $]0, 1]$. Thus by virtue of [4, Corollary 4.3(ii)], it is firmly non-expansive. ☐

Now define $M_{\overline{G}} \stackrel{\text{def}}{=} \mathrm{P}_{T^G_{x^\star}} W_{\overline{G}} \mathrm{P}_{T^G_{x^\star}}$ and $M_{\overline{J}} \stackrel{\text{def}}{=} \mathrm{P}_{T^J_{x^\star}} W_{\overline{J}} \mathrm{P}_{T^J_{x^\star}}$, and the matrices

$$\begin{aligned} M &\stackrel{\text{def}}{=} \mathrm{Id} + 2M_{\overline{G}}M_{\overline{J}} - M_{\overline{G}} - M_{\overline{J}} = \frac{1}{2}\mathrm{Id} + \frac{1}{2}(2M_{\overline{G}} - \mathrm{Id})(2M_{\overline{J}} - \mathrm{Id}), \\ M_\lambda &\stackrel{\text{def}}{=} (1-\lambda)\mathrm{Id} + \lambda M, \quad \lambda \in ]0, 2[. \end{aligned} \tag{5.4}$$

We have the following locally linearized version of (1.5).

**Proposition 5.2 (Locally linearized DR iteration).** *Suppose that the DR iteration (1.5) is run under the assumptions of Theorem 4.3. Assume also that $\lambda_k \to \lambda \in ]0, 2[$. Then $M$ is* firmly non-expansive *and $M_\lambda$ is $\frac{\lambda}{2}$-averaged. Moreover, for all $k$ large enough, we have*

$$z_{k+1} - z^\star = M_\lambda(z_k - z^\star) + \psi_k + \phi_k, \tag{5.5}$$

*where $\|\psi_k\| = o(\|z_k - z^\star\|)$ and $\phi_k = O(\lambda_k|\gamma_k - \gamma|)$. $\psi_k$ and $\phi_k$ vanish when $G$ and $J$ are locally polyhedral around $x^\star$ and $(\gamma_k, \lambda_k)$ are chosen constant in $]0, +\infty[\times]0, 2[$.*

See Section C.2 for the proof.

**Remark 5.3.** If $\phi_k = o(\|z_k - z^\star\|)$, then the rest in (5.5) is $o(\|z_k - z^\star\|)$. However, this is of little practical interest as $z^\star$ is unknown.

We now derive a characterization of the spectral properties of $M_\lambda$, which in turn, will allow to study the linear convergence rates of its powers. Recall the notion of convergent matrices from Definition 2.5. To lighten the notation in the following, we will set $S^J_{x^\star} \stackrel{\text{def}}{=} (T^J_{x^\star})^\perp$ and $S^G_{x^\star} \stackrel{\text{def}}{=} (T^G_{x^\star})^\perp$.

**Lemma 5.4.** *Suppose that $\lambda \in ]0, 2[$, then,*

(i) *$M_\lambda$ is convergent with*

$$M^\infty = \mathrm{P}_{\ker(M_{\overline{G}}(\mathrm{Id}-M_{\overline{J}})+(\mathrm{Id}-M_{\overline{G}})M_{\overline{J}})},$$

*and we have*

$$\forall k \in \mathbb{N},\ M_\lambda^k - M^\infty = (M_\lambda - M^\infty)^k \ \text{ and } \ \rho(M_\lambda - M^\infty) < 1.$$

*In particular, if*

$$T_{x^\star}^J \cap T_{x^\star}^G = \{0\},\ \ \mathrm{span}(\mathrm{Id} - M_{\overline{J}}) \cap S_{x^\star}^G = \{0\} \ \text{ and } \ \mathrm{span}(\mathrm{Id} - M_{\overline{G}}) \cap T_{x^\star}^G = \{0\},$$

*then $M^\infty = 0$.*

(ii) *Given any $\rho \in ]\rho(M_\lambda - M^\infty), 1[$, there is $K$ large enough such that for all $k \geq K$,*

$$\|M_\lambda^k - M^\infty\| = O(\rho^k).$$

(iii) *If, moreover, $G$ and $J$ are locally polyhedral around $x^\star$, then $M_\lambda$ is normal (i.e. $M_\lambda^T M_\lambda = M_\lambda M_\lambda^T$) and converges linearly to $\mathrm{P}_{(T_{x^\star}^J \cap T_{x^\star}^G) \oplus (S_{x^\star}^J \cap S_{x^\star}^G)}$ with the optimal rate*

$$\sqrt{(1-\lambda)^2 + \lambda(2-\lambda)\cos^2\big(\theta_F(T_{x^\star}^J, T_{x^\star}^G)\big)} < 1.$$

*In particular, if $T_{x^\star}^J \cap T_{x^\star}^G = S_{x^\star}^J \cap S_{x^\star}^G = \{0\}$, then $M_\lambda$ converges linearly to $0$ with the optimal rate*

$$\sqrt{(1-\lambda)^2 + \lambda(2-\lambda)\cos^2\big(\theta_1(T_{x^\star}^J, T_{x^\star}^G)\big)} < 1.$$

See Section C.3 for the proof.

Combining Proposition 5.2 and Lemma 5.4, we have the following equivalent characterization of the locally linearized iteration.

**Corollary 5.5.** *Suppose that the DR iteration (1.5) is run under the assumptions of Theorem 4.3. Assume also that $\lambda_k \to \lambda \in ]0, 2[$. Then the following holds.*

(i) *(5.5) is equivalent to*

$$(\mathrm{Id} - M^\infty)(z_{k+1} - z^\star) = (M_\lambda - M^\infty)(\mathrm{Id} - M^\infty)(z_k - z^\star) + (\mathrm{Id} - M^\infty)\psi_k + \phi_k. \tag{5.6}$$

(ii) *If $G$ and $J$ are locally polyhedral around $x^\star$ and $(\gamma_k, \lambda_k)$ are chosen constant in $]0, +\infty[\times]0, 2[$, then*

$$z_{k+1} - z^\star = (M_\lambda - M^\infty)(z_k - z^\star). \tag{5.7}$$

The direction $\Rightarrow$ is easy, the converse needs more arguments. See Section C.4 for the proof.

## 5.2 Local linear convergence

We are now in position to present the local linear convergence of DR method.

**Theorem 5.6 (Local linear convergence of DR).** *Suppose that the DR iteration (1.5) is run under the assumptions of Proposition 5.2. Recall $M^\infty$ from Lemma 5.4. The following holds:*

(i) *Given any $\rho \in ]\rho(M_\lambda - M^\infty), 1[$, there exists $K \in \mathbb{N}$ large enough such that for all $k \geq K$, if $\lambda_k|\gamma_k - \gamma| = O(\eta^k)$ for $0 \leq \eta < \rho$, then*

$$\|(\mathrm{Id} - M^\infty)(z_k - z^\star)\| = O(\rho^{k-K}).$$

(ii) *If $G$ and $J$ are locally polyhedral around $x^\star$ and $(\gamma_k, \lambda_k) \equiv (\gamma, \lambda) \in ]0, +\infty[\times]0, 2[$, then there exists $K \in \mathbb{N}$ large enough such that for all $k \geq K$,*

$$\|z_k - z^\star\| \leq \rho^{k-K} \|z_K - z^\star\| \tag{5.8}$$

*where the convergence rate*

$$\rho = \sqrt{(1-\lambda)^2 + \lambda(2-\lambda)\cos^2\big(\theta_F(T^J_{x^\star}, T^G_{x^\star})\big)} \in [0, 1[$$

*is optimal.*

See Section C.5 for the proof.

**Remark 5.7.**

(i) If $M^\infty = 0$ in (i) or in the situation of (ii), we also have local linear convergence of $x_k$ and $v_k$ to $x^\star$ by non-expansiveness of the proximity operator.

(ii) The condition on $\phi_k$ in Theorem 5.6(i) amounts to saying that $\gamma_k$ should converge fast enough to $\gamma$. Otherwise, the local convergence rate would be dominated by that of $\phi_k$. Especially when $\phi_k$ converges sub-linearly to 0, then the local convergence rate will eventually become sublinear. See Figure 5 in the numerical experiments section.

(iii) For Theorem 5.6(ii), it can be observed that the best rate is obtained for $\lambda = 1$. This has been also pointed out in [19] for basis pursuit. This assertion is however only valid for the local convergence behaviour and does not mean in general that the DR will be globally faster for $\lambda_k \equiv 1$.

(iv) Observe also that the local linear convergence rate does not depend on $\gamma$ when both $G$ and $J$ are locally polyhedral around $x^\star$. This means that the choice of $\gamma_k$ only affects the number of iterations needed for finite identification. For general partly smooth functions, $\gamma_k$ influences both the identification time and the local linear convergence rate, since $M_\lambda$ depends on it through the matrices $W_{\overline{G}}$ and $W_{\overline{J}}$ ($\gamma$ weights the Riemannian Hessians of $\overline{G}$ and $\overline{J}$; see (5.1)-(5.3)). See Figure 4 for a numerical comparison.

# 6 Finite convergence

We are now ready to characterize situations where finite convergence of DR occurs.

**Theorem 6.1.** *Assume that the unrelaxed stationary DR iteration is used (i.e., $\gamma_k \equiv \gamma \in ]0, +\infty[$ and $\lambda_k \equiv 1$), such that $(z_k, x_k, v_k) \to (z^\star, x^\star, x^\star)$, where $G$ and $J$ are locally polyhedral nearby $x^\star$. Suppose that either $J$ or $G$ is locally $C^2$ at $x^\star$. Then the DR sequences $\{z_k, x_k, v_k\}_{k\in\mathbb{N}}$ converge in finitely many steps to $(z^\star, x^\star, x^\star)$.*

**Proof.** We will prove the statement when $J$ is locally $C^2$ at $x^\star$, and the same reasoning holds if the assumption is on $G$. Local $C^2$-smoothness of $J$ at $x^\star$ entails that $\partial J(x^\star) = \{\nabla J(x^\star)\}$ and $J$ is partly smooth at $x^\star$ relative to $\mathcal{M}^J_{x^\star} = \mathbb{R}^n$. Moreover, the non-degeneracy condition (ND) is in force. It then follows from Proposition 5.2 and Lemma 5.4(i) that there exists $K \in \mathbb{N}$ large enough such that

$$\forall k \geq K, \quad z_{k+1} - z^\star = \mathrm{P}_{T^G_{x^\star}}(z_k - z^\star) \Rightarrow \forall k \geq K+1, \quad z_k - z^\star \in T^G_{x^\star},$$

whence we conclude that

$$\forall k \geq K+1, \quad z_k = z_{k+1} = \cdots = z^\star. \qquad \square$$

DR is known (see, *e.g.*, [21, Theorem 6]) to be a special case of the exact proximal point algorithm (PPA) with constant step-size $\gamma_k \equiv 1$. This suggests that many results related to PPA can be carried over to DR. For instance, finite convergence of PPA has been studied in [42, 38] under different conditions.

However, [21, Theorem 9] gave a negative result that suggests that these previous conditions sufficient for finite termination of PPA can be difficult or impossible to carry over to DR even for the polyhedral case. The authors in [5] considered the unrelaxed and stationary DR for solving the convex feasibility problem

$$\text{Find a point in } C_1 \cap C_2,$$

where $C_1$ and $C_2$ are nonempty closed convex sets in $\mathbb{R}^n$, $C_1 \cap C_2 \neq \emptyset$, $C_1$ is an affine subspace and $C_2$ is a polyhedron. They established finite convergence under Slater's condition

$$C_1 \cap \mathrm{int}(C_2) \neq \emptyset.$$

They also provided examples where this condition holds where the conditions of [42, 38] for finite convergence do not apply.

Specializing our result to $G = \iota_{C_1}$ and $J = \iota_{C_2}$, then under Slater's condition, if $x^\star \in C_1 \cap \mathrm{int}(C_2)$, we have $G$ is partly smooth at any $x \in C_1$ relative to $C_1$ with $T^G_{x^\star} = \mathrm{par}(C_1)$ (*i.e.* a translate of $C_1$ to the origin), and $\partial J(x^\star) = N_{C_2}(x^\star) = \{0\}$, and we recover the result of [5]. In fact, [5, Theorem 3.7] shows that the cluster point $x^\star$ is always an interior point regardless of the starting point of DR. The careful reader may have noticed that in the current setting, thanks to Example 4.6, the estimate in (4.5) gives a bound on the finite convergence iteration.

# 7 Alternating direction method of multipliers

For problem ($\mathcal{P}$), let us now compose $J$ with a linear operator $L$, which results in the following problem

$$\min_{x\in\mathbb{R}^n} G(x) + J(Lx), \tag{$\mathcal{P}_L$}$$

where

- **(A.4)** $G \in \Gamma_0(\mathbb{R}^n)$ is the same as problem ($\mathcal{P}$), $J \in \Gamma_0(\mathbb{R}^m)$ and $L : \mathbb{R}^n \to \mathbb{R}^m$ is an *injective* linear operator;
- **(A.5)** $\mathrm{ri}(\mathrm{dom}(G) \cap \mathrm{dom}(J)) \neq \emptyset$;
- **(A.6)** $\mathrm{Argmin}(G + J \circ L) \neq \emptyset$, *i.e.* the set of minimizers is non-empty.

The main difficulty of using DR solve the composed problem ($\mathcal{P}_L$) is that the proximity operator of the composition $J \circ L$ in general can not be solved explicitly. The alternating direction method of multipliers [22] is an efficient to deal with such difficulty. The stationary version (*i.e.* constant step-size) of the method is described in Algorithm 2.

**Algorithm 2:** Alternating Direction method of Multipliers

**Initial**: $\gamma > 0$. $k = 0$, $v_0, y_0 \in \mathbb{R}^m$;
**repeat**

$$\begin{aligned} u_{k+1} &= \mathrm{argmin}_{u\in\mathbb{R}^n} G(u) + \langle Lu,\, y_k\rangle + \frac{\gamma}{2}\|Lu - v_k\|^2, \\ v_{k+1} &= \mathrm{argmin}_{v\in\mathbb{R}^m} J(v) - \langle v,\, y_k\rangle + \frac{\gamma}{2}\|v - Lu_{k+1}\|^2, \\ y_{k+1} &= y_k + \gamma(Lu_{k+1} - v_{k+1}), \end{aligned} \tag{7.1}$$

**until** *convergence*;

It is shown in [22] (see also [21]), that ADMM amounts to applying DR to the Fenchel-Rockafellar dual problem of ($\mathcal{P}_L$). In the following, we recall in short the derivation of transforming ADMM to DR. First, the dual form of ($\mathcal{P}_L$) is

$$\max_{y\in\mathbb{R}^m} -\big(G^*(-L^Ty) + J^*(y)\big), \tag{$\mathcal{D}_L$}$$

where $G^*$, $J^*$ denote the Fenchel conjugate of $G$ and $J$ respectively, and given function $G \in \Gamma_0(\mathbb{R}^n)$, its conjugate is defined as

$$G^*(y) = \sup_{x \in \mathbb{R}^n} \langle y,\, x \rangle - G(x).$$

Define $w_{k+1} = y_k + \gamma L u_{k+1}$, then we get the follow the iteration from (7.1)

$$\begin{cases} u_{k+1} = \mathrm{argmin}_{u \in \mathbb{R}^n} G(u) + \frac{\gamma}{2}\|Lu - \frac{1}{\gamma}(w_k - 2y_k)\|^2, \\ w_{k+1} = y_k + \gamma L u_{k+1}, \\ v_{k+1} = \mathrm{argmin}_{v \in \mathbb{R}^m} J(v) + \frac{\gamma}{2}\|v - \frac{1}{\gamma} w_{k+1}\|^2, \\ y_{k+1} = w_{k+1} - \gamma v_{k+1}. \end{cases} \tag{7.2}$$

Apply the Moreau's identity to $G, J$ respectively [43], then we obtain

$$\begin{cases} w_{k+1} = \mathrm{argmin}_{w \in \mathbb{R}^n} G^*(-L^T w) + \frac{1}{2\gamma}\|w - (2y_k - w_k)\|^2 + (w_k - y_k), \\ y_{k+1} = \mathrm{argmin}_{y \in \mathbb{R}^m} J^*(y) + \frac{1}{2\gamma}\|y - w_{k+1}\|^2, \end{cases} \tag{7.3}$$

which is clearly applying the non-relaxed and stationary DR method (*i.e.* $\lambda_k \equiv 1, \gamma_k \equiv \gamma$) to the dual problem ($\mathcal{D}_L$).

When $L$ injective (*i.e.* has full column rank), the convergence of all the sequences in (7.2) are guaranteed, that is

$$u_k \to u^\star, \;\; v_k \to v^\star = Lu^\star, \;\; w_k \to w^\star, \;\; y_k \to y^\star, \tag{7.4}$$

where $u^\star$ is a global minimizer of ($\mathcal{P}_L$), $v^\star$ is the dual solution of ($\mathcal{D}_L$), $w^\star$ is a fixed point of (7.3), and $y^\star$ is the Lagrangian multiplier.

Owing to the result of [31, Theorem 4.2], the class of partly smooth functions is closed under pre-composition by a surjective linear operator. Hence, if $G^*$ is partly smooth, then so is the composition $G^*(-L^T)$. As a consequence, besides conditions (**A.4**)-(**A.6**), if we assume that $G^*$ and $J^*$ are moreover partly smooth functions, then based on the result of Section 4.1 and 5, we can obtain the local linear convergence of the ADMM through the dual iteration (7.3). As a matter of fact, we can establish result directly to the primal ADMM iteration (7.1), which is the content of the rest of this section.

## 7.1 Local linear convergence of ADMM

Let $u^\star$ be a global minimizer of ($\mathcal{P}_L$), $v^\star$ be a dual solution of ($\mathcal{D}_L$) such that $Lu^\star = v^\star$, and $y^\star$ be the Lagrangian multiplier. Then at convergence, we have the following inclusion from (7.2),

$$-L^T y^\star \in \partial G(u^\star) \quad \text{and} \quad y^\star \in \partial J(v^\star).$$

**Proposition 7.1 (Finite activity identification).** *For the ADMM iteration* (7.1)*, suppose that assumptions* (**A.4**)–(**A.6**) *are hold such that the created sequence* $(u_k, v_k, y_k) \to (u^\star, v^\star, y^\star)$ *where* $u^\star \in \mathrm{argmin}(G + J \circ L)$ *and* $v^\star = Lu^\star$*. If* $G \in \mathrm{PSF}_{u^\star}(\mathcal{M}^G_{u^\star})$, $J \in \mathrm{PSF}_{v^\star}(\mathcal{M}^J_{v^\star})$*, and moreover the non-degeneracy condition*

$$-\,L^T y^\star \in \mathrm{ri}\big(\partial G(u^\star)\big) \;\; \textit{and} \;\; y^\star \in \mathrm{ri}\big(\partial J(v^\star)\big), \tag{ND$_{\mathcal{D}}$}$$

*holds. Then,*

(i) *For all $k$ sufficiently large,* $(u_k, v_k) \in \mathcal{M}^G_{u^\star} \times \mathcal{M}^J_{v^\star}$.

(ii) *Moreover,*

  (a) *If* $\mathcal{M}^G_{u^\star} = u^\star + T^G_{u^\star}$*, then* $\forall k \geq K$, $T^G_{u_k} = T^G_{u^\star}$.

(b) *If $G$ is locally polyhedral around $u^\star$, then $\forall k \geq K$, $u_k \in \mathcal{M}^G_{u^\star} = u^\star + T^G_{u^\star}$, $T^G_{u_k} = T^G_{u^\star}$, $\nabla_{\mathcal{M}^G_{u^\star}} G(u_k) = \nabla_{\mathcal{M}^G_{u^\star}} G(u^\star)$, and $\nabla^2_{\mathcal{M}^G_{u^\star}} G(u_k) = 0$.*

*Similar conclusion holds for function $J$ and sequence $v_k$.*

**Proof.** At convergence, we have $y^\star = w^\star - \gamma L u^\star = w^\star - \gamma v^\star$, $-L^T y^\star = L^T(w^\star - 2y^\star - \gamma L u^\star)$. Then from the update of $u_{k+1}, v_{k+1}$ in (7.2), we have the following monotone inclusions

$$\begin{aligned} L^T(w_k - 2y_k - \gamma L u_{k+1}) &\in \partial G(u_{k+1}) \quad \text{and} \quad w_k - \gamma v_k \in \partial J(v_k), \\ L^T(w^\star - 2y^\star - \gamma L u^\star) &\in \partial G(u^\star) \quad \text{and} \quad w^\star - \gamma v^\star \in \partial J(v^\star), \end{aligned}$$

Since $L$ is bounded linear operator, it then follows that

$$\begin{aligned} \mathrm{dist}\big(-L^T y^\star, \partial G(u_{k+1})\big) &\leq \|L^T(w_k - 2y_k - \gamma L u_{k+1}) - L^T(w^\star - 2y^\star - \gamma L u^\star)\| \\ &\leq \|L\| \|(w_k - w^\star) - 2(y_k - y^\star) - \gamma L(u_{k+1} - u^\star)\| \\ &\leq \|L\| \big(\|w_k - w^\star\| + 2\|y_k - y^\star\| + \gamma \|L\| \|u_{k+1} - u^\star\|\big) \to 0. \end{aligned}$$

and

$$\mathrm{dist}\big(y^\star, \partial J(v_k)\big) \leq \|(w_k - \gamma v_k) - (w^\star - \gamma v^\star)\| \leq \|w_k - w^\star\| + \gamma \|v_k - v^\star\| \to 0.$$

Since $G \in \Gamma_0(\mathbb{R}^n)$ and $J \in \Gamma_0(\mathbb{R}^m)$, then by the sub-differentially continuous property of them we have $G(u_k) \to G(u^\star)$ and $J(v_k) \to J(v^\star)$. Hence the conditions of [24, Theorem 5.3] are fulfilled for $G$ and $J$, and the finite identification claim follows. The rest of the proof follows the proof of Theorem 4.3. □

Define the following function similar to those in (5.1),

$$\overline{\overline{G}}(u) \stackrel{\text{def}}{=} \frac{1}{\gamma}\big(G(u) - \langle u,\, -L^T y^\star \rangle\big), \quad \overline{\overline{J}}(v) \stackrel{\text{def}}{=} \frac{1}{\gamma}\big(J(v) - \langle v,\, y^\star \rangle\big). \tag{7.5}$$

Owing to Lemma 4.2, with condition ($\mathrm{ND}_{\mathcal{D}}$) holding, their Riemannian Hessian are positive semidefinite. Hence, define the following matrices

$$\begin{aligned} H_{\overline{\overline{G}}} &\stackrel{\text{def}}{=} P_{T^G_{u^\star}} \nabla^2_{\mathcal{M}^G_{u^\star}} \overline{\overline{G}}(u^\star) P_{T^G_{u^\star}} \quad &&\text{and} \quad H_{\overline{\overline{J}}} \stackrel{\text{def}}{=} P_{T^J_{v^\star}} \nabla^2_{\mathcal{M}^J_{v^\star}} \overline{\overline{J}}(v^\star) P_{T^J_{v^\star}}, \\ W_{\overline{\overline{G}}} &\stackrel{\text{def}}{=} L_G\big(\mathrm{Id} + (L_G^T L_G)^{-1} H_{\overline{\overline{G}}}\big)^{-1} (L_G^T L_G)^{-1} L_G^T \quad &&\text{and} \quad W_{\overline{\overline{J}}} \stackrel{\text{def}}{=} \big(\mathrm{Id} + H_{\overline{\overline{J}}}\big)^{-1}, \end{aligned} \tag{7.6}$$

where $L_G \stackrel{\text{def}}{=} L \circ P_{T^G_{u^\star}}$. Define $T^{G,L}_{u^\star}$ the subspace corresponding to the projection operator $P_{T^{G,L}_{u^\star}} = L_G(L_G^T L_G)^{-1} L_G^T$.

Now define $M_{\overline{\overline{G}}} \stackrel{\text{def}}{=} P_{T^G_{u^\star}} W_{\overline{\overline{G}}} P_{T^G_{u^\star}}$ and $M_{\overline{\overline{J}}} \stackrel{\text{def}}{=} P_{T^J_{v^\star}} W_{\overline{\overline{J}}} P_{T^J_{v^\star}}$, and the matrix

$$M = \mathrm{Id} + 2M_{\overline{\overline{G}}} M_{\overline{\overline{J}}} - M_{\overline{\overline{G}}} - M_{\overline{\overline{J}}} = \frac{1}{2}(2M_{\overline{\overline{G}}} - \mathrm{Id})(2M_{\overline{\overline{J}}} - \mathrm{Id}) + \frac{1}{2}\mathrm{Id}. \tag{7.7}$$

From Lemma 5.4, we have $M$ is convergent.

**Proposition 7.2.** *Let (A.4)-(A.6) and conditions in Theorem 7.1 hold.*

(i) *Matrix $M$ is firmly non-expansive, hence convergent.*

(ii) *For all $k$ large enough, the ADMM iteration (7.1) can be written as the following fixed-point iteration*

$$w_{k+1} - w^\star = M(w_k - w^\star) + o(\|w_k - w^\star\|). \tag{7.8}$$

(iii) *For the sequence $\{w_k\}_{k\in\mathbb{N}}$, we have the following result*

(a) *given any $\rho \in ]\rho(M - M^\infty), 1[$, there is $K$ large enough such that for all $k \geq K$,*

$$\|(\mathrm{Id} - M^\infty)(w_k - w^\star)\| = O(\rho^{k-K}). \tag{7.9}$$

(b) *If $G$ and $J$ are locally polyhedral around $u^\star$ and $v^\star$ respectively, then there exists $K \in \mathbb{N}$ large enough such that for all $k \geq K$,*

$$\|w_k - w^\star\| \leq \rho^{k-K}\|w_K - w^\star\|, \tag{7.10}$$

*where the value of $\rho$ is*

$$\rho = \cos\big(\theta_F(T^J_{v^\star}, T^{G,L}_{u^\star})\big) \in [0, 1[,$$

(iv) *Moreover, if $M^\infty = 0$, given any $\rho \in ]\rho(M), 1[$, there is $K$ large enough such that for all $k \geq K$,*

$$\|v_k - v^\star\| = O(\rho^k),\ \|y_k - y^\star\| = O(\rho^k),\ \|L(u_{k+1} - u^\star)\| = O(\rho^k). \tag{7.11}$$

*If moreover $G, J$ are locally polyhedral around $u^\star, v^\star$ respectively, then (7.11) holds with*

$$\rho = \cos\big(\theta_F(T^J_{v^\star}, T^{G,L}_{u^\star})\big) \in [0, 1[,$$

*which is the optimal convergence rate.*

See Section C.6 for the proof of the above proposition.

## 8 Sum of more than two functions

We now want to tackle the problem of solving

$$\min_{x\in\mathbb{R}^n} \textstyle\sum_{i=1}^m J_i(x), \tag{$\mathcal{P}_m$}$$

where

**(A'.1)** $J_i \in \Gamma_0(\mathbb{R}^n)$, $\forall i = 1, \cdots, m$;

**(A'.2)** $\bigcap_{1\leq i\leq m} \mathrm{ri}(\mathrm{dom}(J_i)) \neq \emptyset$;

**(A'.3)** $\mathrm{Argmin}(\sum_{i=1}^m J_i) \neq \emptyset$.

In fact, problem ($\mathcal{P}_m$) can be equivalently reformulated as ($\mathcal{P}$) in a product space, see *e.g.* [14, 41]. Let $\boldsymbol{\mathcal{H}} = \underbrace{\mathbb{R}^n \times \cdots \times \mathbb{R}^n}_{m \text{ times}}$ endowed with the scalar inner-product and norm

$$\forall \boldsymbol{x}, \boldsymbol{y} \in \boldsymbol{\mathcal{H}},\ \langle \boldsymbol{x}, \boldsymbol{y}\rangle = \textstyle\sum_{i=1}^m \langle x_i,\, y_i\rangle,\ \|\boldsymbol{x}\| = \sqrt{\sum_{i=1}^m \|x_i\|^2}.$$

Let $\boldsymbol{\mathcal{S}} = \big\{\boldsymbol{x} = (x_i)_i \in \boldsymbol{\mathcal{H}} :\ x_1 = \cdots = x_m\big\}$ and its orthogonal complement $\boldsymbol{\mathcal{S}}^\perp = \{\boldsymbol{x} = (x_i)_i \in \boldsymbol{\mathcal{H}} : \sum_{i=1}^m x_i = 0\}$. Now define the canonical isometry,

$$\boldsymbol{C} : \mathbb{R}^n \to \boldsymbol{\mathcal{S}},\ x \mapsto (x, \cdots, x),$$

then we have $\mathrm{P}_{\boldsymbol{\mathcal{S}}}(\boldsymbol{z}) = \boldsymbol{C}(\frac{1}{m}\sum_{i=1}^m z_i)$.

Problem ($\mathcal{P}_m$) is now equivalent to

$$\min_{\boldsymbol{x}\in\boldsymbol{\mathcal{H}}} \boldsymbol{J}(\boldsymbol{x}) + \boldsymbol{G}(\boldsymbol{x}),\ \text{ where } \boldsymbol{J}(\boldsymbol{x}) = \textstyle\sum_{i=1}^m J_i(x_i) \text{ and } \boldsymbol{G}(\boldsymbol{x}) = \iota_{\boldsymbol{\mathcal{S}}}(\boldsymbol{x}), \tag{$\boldsymbol{\mathcal{P}}$}$$

which has the same structure on $\boldsymbol{\mathcal{H}}$ as ($\mathcal{P}$) on $\mathbb{R}^n$.

Obviously, $\boldsymbol{J}$ is separable and therefore,

$$\mathrm{prox}_{\gamma\boldsymbol{J}}(\boldsymbol{x}) = \big(\mathrm{prox}_{\gamma J_i}(x_i)\big)_i.$$

Let $\boldsymbol{x}^\star = \boldsymbol{C}(x^\star)$. Clearly, $\boldsymbol{G}$ is polyhedral, hence partly smooth relative to $\boldsymbol{\mathcal{S}}$ with $\boldsymbol{T}_{\boldsymbol{x}^\star}^{\boldsymbol{G}} = \boldsymbol{\mathcal{S}}$. Suppose that $J_i \in \mathrm{PSF}_{x^\star}(\mathcal{M}_{x^\star}^{J_i})$ for each $i$. Denote $\boldsymbol{T}_{\boldsymbol{x}^\star}^{\boldsymbol{J}} = \times_i T_{x^\star}^{J_i}$ and $\boldsymbol{S}_{\boldsymbol{x}^\star}^{\boldsymbol{J}} = \times_i (T_{x^\star}^{J_i})^\perp$. Similarly to (5.2), define

$$\boldsymbol{Q} \stackrel{\text{def}}{=} \mathrm{P}_{\boldsymbol{T}_{\boldsymbol{x}^\star}^{\boldsymbol{J}}} \nabla^2 \overline{\boldsymbol{J}}(x^\star) \mathrm{P}_{\boldsymbol{T}_{\boldsymbol{x}^\star}^{\boldsymbol{J}}} \quad \text{and} \quad \boldsymbol{V} \stackrel{\text{def}}{=} (\mathbf{Id} + \boldsymbol{Q})^{-1},$$

where $\overline{\boldsymbol{J}}(\boldsymbol{x}) \stackrel{\text{def}}{=} \gamma \sum_{i=1}^m \widetilde{J}_i(x_i) - \langle \boldsymbol{x}, \boldsymbol{z}^\star - \boldsymbol{x}^\star \rangle$, $\widetilde{J}_i$ is the smooth representation of $J_i$ on $\mathcal{M}_{x^\star}^{J_i}$, and $\mathbf{Id}$ is the identity operatror on $\boldsymbol{\mathcal{H}}$. Since $\boldsymbol{J}$ is polyhedeal, we have $\boldsymbol{U} = \mathbf{Id}$. Now we can provide the product space form of (7.7), which reads

$$\begin{aligned} \boldsymbol{M} &= \mathbf{Id} + 2\mathrm{P}_{\boldsymbol{T}_{\boldsymbol{x}^\star}^{\boldsymbol{G}}} \mathrm{P}_{\boldsymbol{T}_{\boldsymbol{x}^\star}^{\boldsymbol{J}}} \boldsymbol{V} \mathrm{P}_{\boldsymbol{T}_{\boldsymbol{x}^\star}^{\boldsymbol{J}}} - \mathrm{P}_{\boldsymbol{T}_{\boldsymbol{x}^\star}^{\boldsymbol{G}}} - \mathrm{P}_{\boldsymbol{T}_{\boldsymbol{x}^\star}^{\boldsymbol{J}}} \boldsymbol{V} \mathrm{P}_{\boldsymbol{T}_{\boldsymbol{x}^\star}^{\boldsymbol{J}}} \\ &= \frac{1}{2}\mathbf{Id} + \mathrm{P}_{\boldsymbol{T}_{\boldsymbol{x}^\star}^{\boldsymbol{G}}} (2\mathrm{P}_{\boldsymbol{T}_{\boldsymbol{x}^\star}^{\boldsymbol{J}}} \boldsymbol{V} \mathrm{P}_{\boldsymbol{T}_{\boldsymbol{x}^\star}^{\boldsymbol{J}}} - \mathbf{Id}) - \frac{1}{2}(2\mathrm{P}_{\boldsymbol{T}_{\boldsymbol{x}^\star}^{\boldsymbol{J}}} \boldsymbol{V} \mathrm{P}_{\boldsymbol{T}_{\boldsymbol{x}^\star}^{\boldsymbol{J}}} - \mathbf{Id}) \\ &= \frac{1}{2}\mathbf{Id} + \frac{1}{2}(2\mathrm{P}_{\boldsymbol{T}_{\boldsymbol{x}^\star}^{\boldsymbol{G}}} - \mathbf{Id})(2\mathrm{P}_{\boldsymbol{T}_{\boldsymbol{x}^\star}^{\boldsymbol{J}}} \boldsymbol{V} \mathrm{P}_{\boldsymbol{T}_{\boldsymbol{x}^\star}^{\boldsymbol{J}}} - \mathbf{Id}), \end{aligned} \tag{8.1}$$

and $\boldsymbol{M}_\lambda \stackrel{\text{def}}{=} (1-\lambda)\mathbf{Id} + \lambda \boldsymbol{M}$. Owing to Lemma 5.4, we have

$$\boldsymbol{M}^\infty = \mathrm{P}_{\ker(\mathrm{P}_{\boldsymbol{T}_{\boldsymbol{x}^\star}^{\boldsymbol{G}}}(\mathbf{Id} - \mathrm{P}_{\boldsymbol{T}_{\boldsymbol{x}^\star}^{\boldsymbol{J}}} \boldsymbol{V} \mathrm{P}_{\boldsymbol{T}_{\boldsymbol{x}^\star}^{\boldsymbol{J}}}) + (\mathbf{Id} - \mathrm{P}_{\boldsymbol{T}_{\boldsymbol{x}^\star}^{\boldsymbol{G}}})\mathrm{P}_{\boldsymbol{T}_{\boldsymbol{x}^\star}^{\boldsymbol{J}}} \boldsymbol{V} \mathrm{P}_{\boldsymbol{T}_{\boldsymbol{x}^\star}^{\boldsymbol{J}}})},$$

and when all $J_i$'s are locally polyhedral nearby $x^\star$, $\boldsymbol{M}^\infty$ specializes to

$$\boldsymbol{M}^\infty = \mathrm{P}_{(\boldsymbol{T}_{x^\star}^{\boldsymbol{J}} \cap \boldsymbol{\mathcal{S}}) \oplus (\boldsymbol{S}_{x^\star}^{\boldsymbol{J}} \cap \boldsymbol{\mathcal{S}}^\perp)}.$$

**Corollary 8.1.** *Suppose that (A'.1)-(A'.3) and (H.2)-(H.3) holds. Consider the sequence* $\{\boldsymbol{z}_k, \boldsymbol{x}_k, \boldsymbol{v}_k\}_{k\in\mathbb{N}}$ *provided by the nonstationary DR method (1.5) applied to solve ($\boldsymbol{\mathcal{P}}$). Then,*

(i) $(\boldsymbol{z}_k, \boldsymbol{x}_k, \boldsymbol{v}_k)$ *converges to* $(\boldsymbol{z}^\star, \boldsymbol{x}^\star, \boldsymbol{x}^\star)$, *where* $\boldsymbol{x}^\star = \boldsymbol{C}(x^\star)$ *and* $x^\star$ *is a global minimizer of* ($\mathcal{P}_m$).

(ii) *Assume, moreover, that* $\gamma_k \geq \underline{\gamma} > 0$ *and* $\gamma_k \to \gamma$, $J_i \in \mathrm{PSF}_{x^\star}(\mathcal{M}_{x^\star}^{J_i})$ *and*

$$\boldsymbol{z}^\star \in \boldsymbol{x}^\star + \gamma \mathrm{ri}\left(\times_i \partial J_i(x^\star)\right) \cap \boldsymbol{\mathcal{S}}^\perp. \tag{ND}$$

*Then,*

(a) *for all $k$ large enough,* $\boldsymbol{x}_k \in \times_i \mathcal{M}_{x^\star}^{J_i}$.

(b) *in addition, if* $\lambda_k \to \lambda \in ]0,2[$, *then given any* $\rho \in ]\rho(\boldsymbol{M}_\lambda - \boldsymbol{M}^\infty), 1[$, *there exists* $K \in \mathbb{N}$ *large enough such that for all* $k \geq K$, *if* $\lambda_k|\gamma_k - \gamma| = O(\eta^k)$ *where* $0 \leq \eta < \rho$, *then*

$$\|(\mathbf{Id} - \boldsymbol{M}^\infty)(\boldsymbol{z}_k - \boldsymbol{z}^\star)\| = O(\rho^{k-K}).$$

*In particular, if all $J_i$'s are locally polyhedral around $x^\star$ and* $(\gamma_k, \lambda_k) \equiv (\gamma, \lambda) \in ]0, +\infty[\times]0,2[$, *then* $\boldsymbol{z}_k$ *(resp.* $x_k \stackrel{\text{def}}{=} \frac{1}{m}\sum_{i=1}^m \boldsymbol{x}_{k,i}$*) converges locally linearly to* $\boldsymbol{z}^\star$ *(resp.* $x^\star$*) at the optimal rate* $\rho = \sqrt{(1-\lambda)^2 + \lambda(2-\lambda)\cos^2(\theta_F(\boldsymbol{T}_{x^\star}^{\boldsymbol{J}}, \boldsymbol{\mathcal{S}}))} \in [0,1[$.

**Proof.**

(i) Apply Theorem 3.1 to ($\boldsymbol{\mathcal{P}}$).

(ii) (a) By the separability rule, we have $\boldsymbol{J} \in \mathrm{PSF}_{\boldsymbol{x}^\star}(\times_i \mathcal{M}_{x^\star}^{J_i})$, see [31, Proposition 4.5]. We have also $\partial \boldsymbol{G}(\boldsymbol{x}^\star) = N_{\boldsymbol{\mathcal{S}}}(\boldsymbol{x}^\star) = \boldsymbol{\mathcal{S}}^\perp$. Then (**ND**) is simply a specialization of condition (ND) to problem ($\boldsymbol{\mathcal{P}}$). The claim then follows from Theorem 4.3.

(b) This is a direct consequence of Theorem 5.6. For the local linear convergence of $x_k$ to $x^\star$ in the last part, observe that

$$\begin{aligned} \|x_k - x^\star\|^2 &= \|\tfrac{1}{m}\textstyle\sum_{i=1}^m \boldsymbol{x}_{k,i} - \tfrac{1}{m}\sum_{i=1}^m \boldsymbol{x}_i^\star\|^2 \\ &\leq \tfrac{1}{m}\textstyle\sum_{i=1}^m \|\boldsymbol{x}_{k,i} - \boldsymbol{x}_i^\star\|^2 \\ &= \tfrac{1}{m}\textstyle\sum_{i=1}^m \|\mathrm{prox}_{\gamma J_i}(\boldsymbol{z}_{k,i}) - \mathrm{prox}_{\gamma J_i}(\boldsymbol{z}_i^\star)\|^2 \\ &\leq \tfrac{1}{m}\textstyle\sum_{i=1}^m \|\boldsymbol{z}_{k,i} - \boldsymbol{z}_i^\star\|^2 = \tfrac{1}{m}\|\boldsymbol{z}_k - \boldsymbol{z}^\star\|^2. \end{aligned}$$ □

We also have the following corollary of Theorem 6.1.

**Corollary 8.2.** *Assume that the unrelaxed stationary DR iteration is used (i.e.,* $\gamma_k \equiv \gamma \in ]0, +\infty[$ *and* $\lambda_k \equiv 1$*), such that* $(\boldsymbol{z}_k, \boldsymbol{x}_k, \boldsymbol{v}_k) \to (\boldsymbol{z}^\star, \boldsymbol{C}(x^\star), \boldsymbol{C}(x^\star))$*, where,* $\forall i$*,* $J_i$ *is locally polyhedral nearby* $x^\star$ *and is differentiable at* $x^\star$*. Then the sequences* $\{\boldsymbol{z}_k, \boldsymbol{x}_k, \boldsymbol{v}_k, \frac{1}{m}\sum_{i=1}^m \boldsymbol{x}_{k,i}\}_{k\in\mathbb{N}}$ *converge in finitely many steps to* $(\boldsymbol{z}^\star, \boldsymbol{C}(x^\star), \boldsymbol{C}(x^\star), x^\star)$*.*

# 9 Numerical experiments

## 9.1 Examples of tested partly smooth functions

Table 1 provides some examples of partly smooth functions that we will use throughout this section in our numerical experiments. These functions are widely used in the literature to regularize a variety of problems in signal/image processing, machine learning and statistics, see *e.g.* [45] and references therein for details. The corresponding Riemannian gradients can also be found in [45]. Since the $\ell_1$, $\ell_\infty$ and the (anisotropic) total variation semi-norm are polyhedral, their Riemannian Hessian vanishes. The Riemannian Hessians for the $\ell_{1,2}$ and the nuclear norm are also provided in [45].

Table 1: Examples of partly smooth functions. $D_{\mathrm{DIF}}$ stands for the finite differences operator.

| Function | Expression | Partial smooth manifold |
|---|---|---|
| $\ell_1$-norm | $\|x\|_1 = \sum_{i=1}^n \|x_i\|$ | $\mathcal{M} = T_x = \{z \in \mathbb{R}^n : I_z \subseteq I_x\}$, $I_x = \{i : x_i \neq 0\}$ |
| $\ell_{1,2}$-norm | $\sum_{i=1}^m \|x_{b_i}\|$ | $\mathcal{M} = T_x = \{z \in \mathbb{R}^n : I_z \subseteq I_x\}$, $I_x = \{i : x_{b_i} \neq 0\}$ |
| $\ell_\infty$-norm | $\max_{i=1,\ldots,n} \|x_i\|$ | $\mathcal{M} = T_x = \{z \in \mathbb{R}^n : z_{I_x} \in \mathbb{R}\,\mathrm{sign}(x_{I_x})\}$, $I_x = \{i : \|x_i\| = \|x\|_\infty\}$ |
| TV semi-norm | $\|x\|_{\mathrm{TV}} = \|D_{\mathrm{DIF}} x\|_1$ | $\mathcal{M} = T_x = \{z \in \mathbb{R}^n : I_{D_{\mathrm{DIF}} z} \subseteq I_{D_{\mathrm{DIF}} x}\}$, $I_{D_{\mathrm{DIF}} x} = \{i : (D_{\mathrm{DIF}} x)_i \neq 0\}$ |
| Nuclear norm | $\|x\|_* = \sum_{i=1}^r \sigma(x)$ | $\mathcal{M} = \{z \in \mathbb{R}^{n_1 \times n_2} : \mathrm{rank}(z) = \mathrm{rank}(x) = r\}$, $\sigma(x)$ singular values of $x$ |

**Affinely-constrained minimization** Let us first consider the affine-constrained minimization problem

$$\min_{x\in\mathbb{R}^n} J(x) \quad \text{subject to} \quad Lx_{\mathrm{ob}} = Lx, \tag{9.1}$$

where $L : \mathbb{R}^n \to \mathbb{R}^m$ is a linear operator, $x_{\mathrm{ob}} \in \mathbb{R}^n$ is known and $J \in \Gamma_0(\mathbb{R}^n)$. Problem (9.1) is of importance in various areas to find regularized solutions to linear equations (one can think for instance of the active area of compressed sensing, matrix completion, and so on). By identifying $G$ with the indicator function of the affine constraint $C \stackrel{\text{def}}{=} \{x \in \mathbb{R}^n : Lx_{\mathrm{ob}} = Lx\} = x_{\mathrm{ob}} + \ker(L)$, it is immediate to see that $G$ is indeed polyhedral and partly smooth at any $x \in C$ relative to $C$.

We here solve (9.1) with $J$ being the $\ell_1$, $\ell_\infty$, $\ell_{1,2}$ and nuclear norms. For all these cases, the proximity operator of $J$ can be computed very easily. In all these experiments, $L$ is drawn randomly from the standard Gaussian ensemble, *i.e.* compressed sensing/matrix completion scenario, with the following settings:

- **$\ell_1$-norm** $(m, n) = (48, 128)$, $x_{\mathrm{ob}}$ is sparse with 8 nonzero entries;
- **$\ell_{1,2}$-norm** $(m, n) = (48, 128)$, $x_{\mathrm{ob}}$ has 3 nonzero blocks of size 4;
- **$\ell_\infty$-norm** $(m, n) = (123, 128)$, $x_{\mathrm{ob}}$ has 10 saturating components;
- **Nuclear norm** $(m, n) = (500, 1024)$, $x_{\mathrm{ob}} \in \mathbb{R}^{32\times 32}$ and $\mathrm{rank}(x_{\mathrm{ob}}) = 4$.

For each setting, the number of measurements is sufficiently large so that one can prove that the minimizer $x^\star$ is unique, and in particular that $\ker(L) \cap T^J_{x^\star} = \{0\}$ (with high probability); see *e.g.* [46]. We also checked that $\ker(L)^\perp \cap S^J_{x^\star} = \{0\}$, which is equivalent to the uniqueness of the fixed point and also implies that $M^\infty = 0$ (see Lemma 5.4(i)). Thus (ND) is fulfilled, and Theorem 5.6 applies. DR is run in its stationary version (*i.e.* constant $\gamma$).

Figure 1 displays the profile of $\|z_k - z^\star\|$ as a function of $k$, and the starting point of the dashed line is the iteration number at which the active partial smoothness manifold of $J$ is identified (recall that $\mathcal{M}^G_{x^\star} = C$ which is trivially identified from the first iteration). One can easily see that for the $\ell_1$ and $\ell_\infty$ norms, the observed linear convergence coincides with the optimal rate predicted by Theorem 5.6(ii). For the case of $\ell_{1,2}$-norm and nuclear norm, though not optimal, our estimates are very tight.

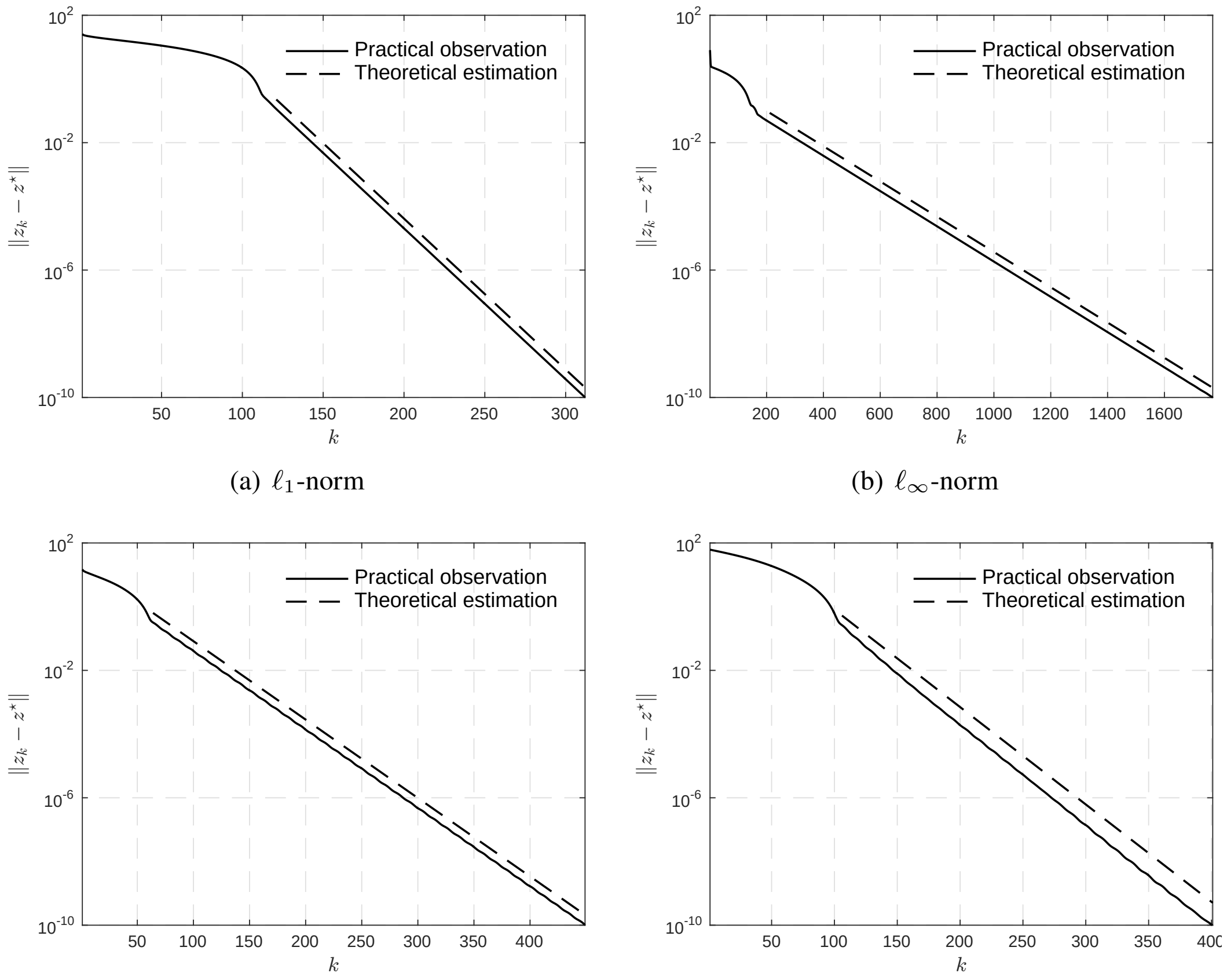

(a) $\ell_1$-norm (b) $\ell_\infty$-norm

(c) $\ell_{1,2}$-norm (d) Nuclear norm

Figure 1: Observed (solid) and predicted (dashed) convergence profiles of DR (1.1) in terms of $\|z_k - z^\star\|$. (a) $\ell_1$-norm. (b) $\ell_\infty$-norm. (c) $\ell_{1,2}$-norm. (d) Nuclear norm. The starting point of the dashed line is the iteration at which the active manifold of $J$ is identified.

**Noise removal** In the following two examples, we suppose that we observe $y = x_{\text{ob}} + \varepsilon$, where $x_{\text{ob}}$ is a piecewise-constant vector, and $\varepsilon$ is an unknown noise supposed to be either uniform or sparse. The goal is to recover $x_{\text{ob}}$ from $y$ using the prior information on $x_{\text{ob}}$ (*i.e.* piecewise-smooth) and $\varepsilon$ (uniform or sparse). To achieve this goal, a popular and natural approach in the signal processing literature is to solve

$$\min_{x \in \mathbb{R}^n} \|x\|_{\text{TV}} \quad \text{subject to} \quad \|y - x\|_p \le \tau, \tag{9.2}$$

where $p = +\infty$ for uniform noise, and $p = 1$ for sparse noise, and $\tau > 0$ is a parameter to be set by the user to adapt to the noise level. Identifying $J = \|\cdot\|_{\text{TV}}$ and $G = \iota_{\|y - \cdot\|_p \le \tau}$, one recognises that for $p \in \{1, +\infty\}$, $J$ and $G$ are indeed polyhedral and their proximity operators are simple to compute.

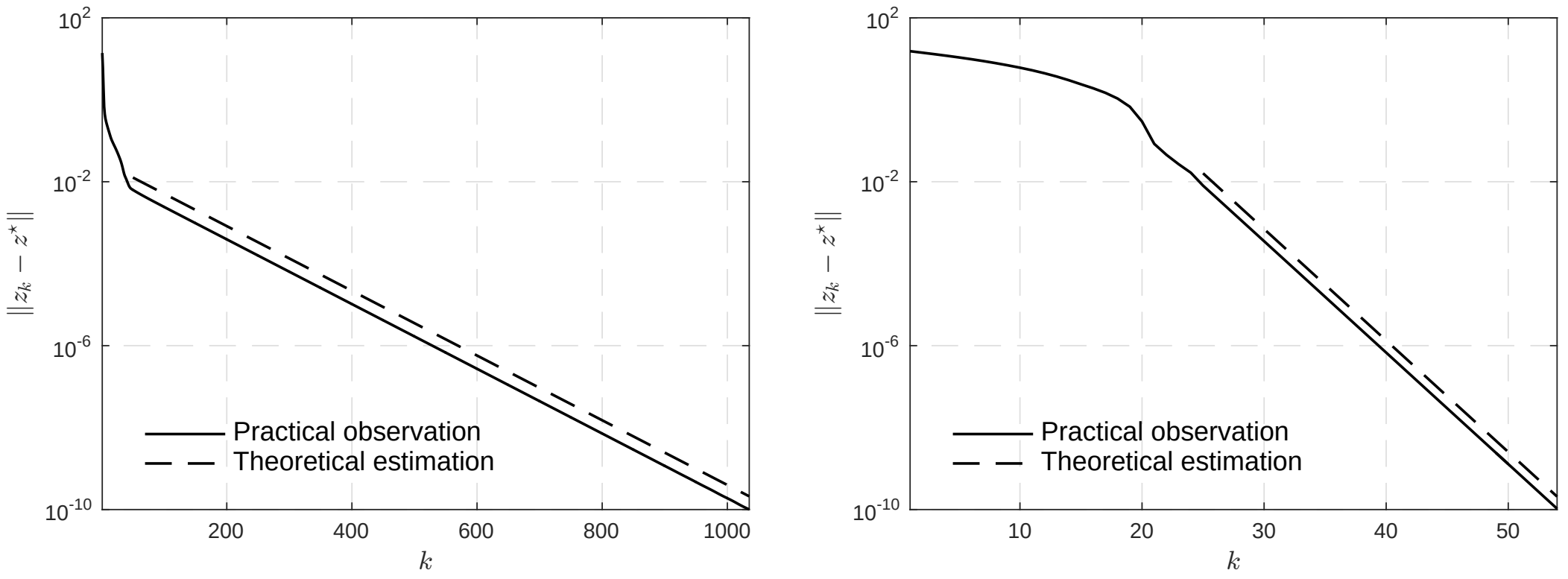

(a) Uniform noise removal $p = +\infty$ (b) Sparse noise removal $p = 1$

Figure 2: Observed (solid) and predicted (dashed) convergence profiles of DR (1.1) in terms of $\|z_k - z^\star\|$. (a) Uniform noise removal by solving (9.2) with $p = +\infty$, (c) Sparse noise removal by solving (9.2) with $p = 1$. The starting point of the dashed line is the iteration at which the manifolds $\mathcal{M}^J_{x^\star}$ and $\mathcal{M}^G_{x^\star}$ are identified.

For both examples, we set $n = 128$ and $x_{\mathrm{ob}}$ is such that $D_{\mathrm{DIF}} x_{\mathrm{ob}}$ has 8 nonzero entries. For $p = +\infty$, $\varepsilon$ is generated uniformly in $[-1, 1]$, and for $p = 1$ $\varepsilon$ is sparse with 16 nonzero entries. DR is run in its stationary version. The corresponding local convergence profiles are depicted in Figure 2(a)-(b). Condition (ND) is checked posterior, and it is satisfied for the considered examples. Owing to polyhedrality, our rate predictions are again optimal.

## 9.2 Finite convergence

We now numerically illustrate the finite convergence of DR. For the remainder of this subsection, we set $n = 2$, and solve ($\mathcal{P}$) with $G = \|\cdot\|_1$ and $J = \iota_C$, $C = \big\{x \in \mathbb{R}^2 : \ \|x - (3/4 \quad 3/4)^T\|_1 \leq 1/2\big\}$. The set of minimizers is the segment $[(1/4 \quad 3/4)^T, (3/4 \quad 1/4)^T]$, and $G$ is differentiable at any minimizer with gradient $(1 \quad 1)^T$. The set of fixed points is thus $[(1/4 \quad 3/4)^T, (3/4 \quad 1/4)^T] - \gamma$. Figure 3(a) shows the trajectory of the sequence $\{z_k\}_{k\in\mathbb{N}}$ and the shadow sequence $\{x_k\}_{k\in\mathbb{N}}$ which both converge finitely as predicted by Theorem 6.1 (DR is used with $\gamma = 0.25$).

For each starting point $z_0 \in [-10, 10]^2$, we run the DR algorithm until $z_{k+1} = z_k$ (up to machine precision), with $\gamma = 0.25$ and $\gamma = 5$. Figure 3(b)-(c) show the number of iterations to finite convergence, where $\gamma = 0.25$ for (b) and $\gamma = 5$ for (c). This confirms that DR indeed converges in finitely many iterations regardless of the starting point and choice of $\gamma$, though more iterations are needed for higher $\gamma$ in this example (see next subsection for further discussion on the choice of $\gamma$).

## 9.3 Choice of $\gamma$

**Impact of $\gamma$ on identification** We now turn to the impact of the choice of $\gamma$ in the DR algorithm. We consider (9.1) with $J$ being the $\ell_1$, the $\ell_{1,2}$ and nuclear norms.

The results are shown in Figure 4, where $K$ denotes the number of iterations needed to identify $\mathcal{M}^J_{x^\star}$ and $\rho$ denotes the local linear convergence rate. We summarize our observations as follows:

- For all examples, the choice of $\gamma$ affects the iteration $K$ at which activity identification occurs. Indeed, $K$ typically decreases monotonically and then either stabilizes or slightly increases. This is in agreement with the bound in (4.1);
- When $J$ is the $\ell_1$, which is polyhedral, the local linear convergence rate is insensitive to $\gamma$ as anticipated by Theorem 5.6(ii). For the other two norms, the local rate depends on $\gamma$ (see Theorem 5.6(i)),

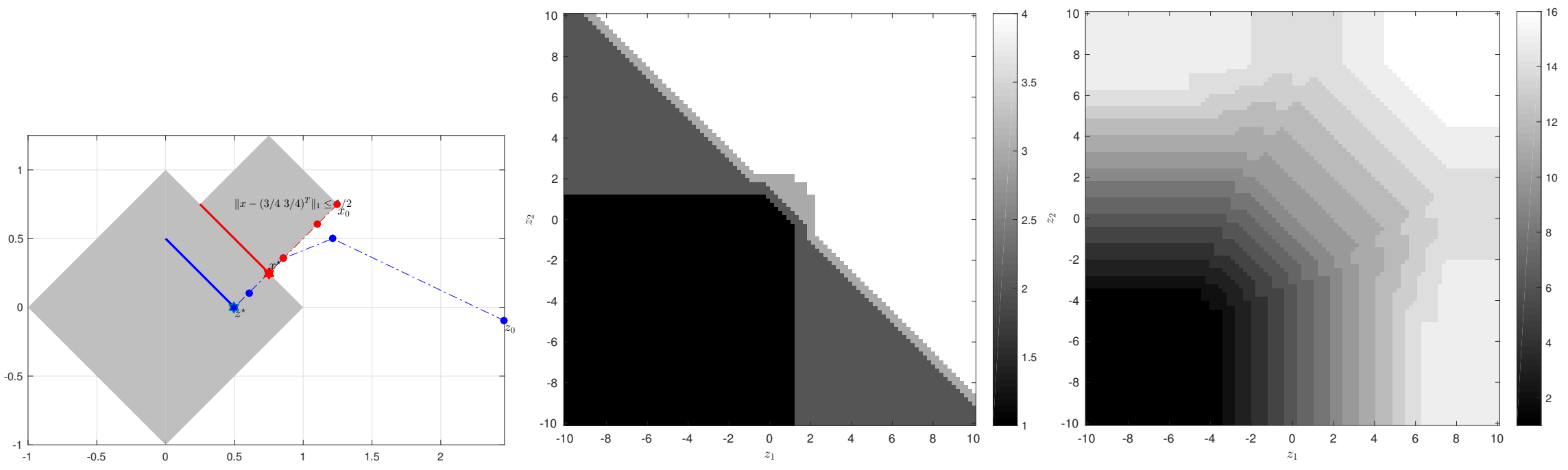

(a) (b) $\gamma = 0.25$ (c) $\gamma = 5$

Figure 3: (a) Trajectories of $\{z_k\}_{k\in\mathbb{N}}$ and $\{x_k\}_{k\in\mathbb{N}}$. The red segment is the set of minimizers and the blue one is the set of fixed points. (b)-(c) Number of iterations needed for the finite convergence of $z_k$ to $z^\star$. DR is run with $\gamma = 0.25$ for (b) and $\gamma = 5$ for (c).

and this rate can be optimized for the parameter $\gamma$;

- In general, there is no correspondence between the optimal choice of $\gamma$ for identification and the one for local convergence rate.

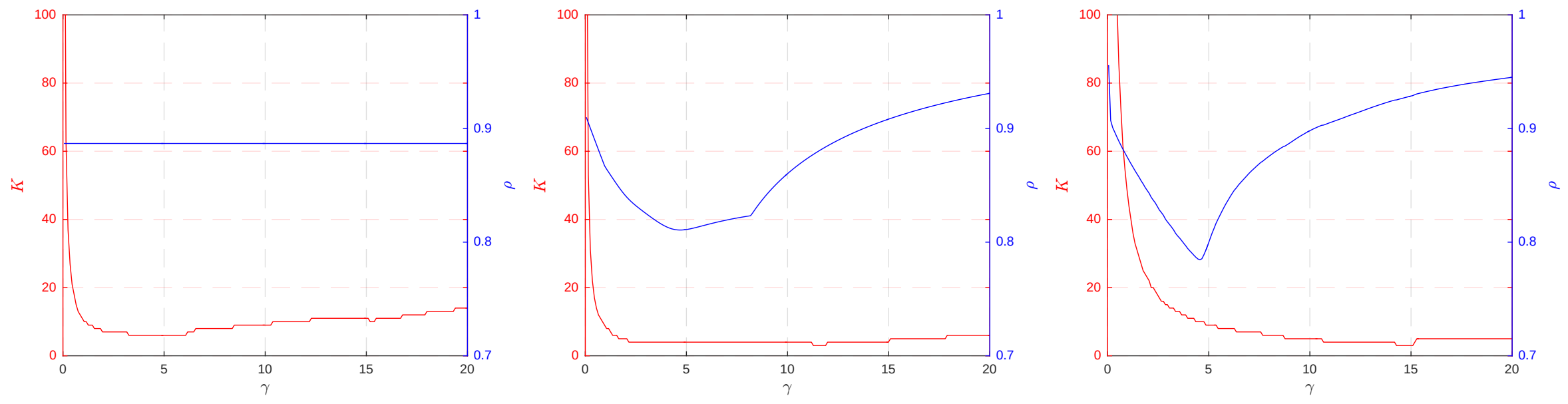

(a) $\ell_1$-norm (b) $\ell_{1,2}$-norm (c) Nuclear norm

Figure 4: Number of iterations ($K$) needed for identification and local linear convergence rate ($\rho$) as a function of $\gamma$ when solving problem (9.1) with different functions $J$ in Table 1. (a) $\ell_1$-norm. (b) $\ell_{1,2}$-norm. (c) Nuclear norm.

**Stationary vs non-stationary DR** We now investigate numerically the convergence behaviour of the non-stationary version of DR and compare it to the stationary one. We fix $\lambda_k \equiv 1$, *i.e.* the iteration is unrelaxed. The stationary DR algorithm is run with some $\gamma > 0$. For the non-stationary one, four choices of $\gamma_k$ are considered:

$$\text{Case 1: } \gamma_k = \gamma + \frac{1}{k^{1.1}}, \ \text{Case 2: } \gamma_k = \gamma + \frac{1}{k^2}, \ \text{Case 3: } \gamma_k = \gamma + 0.95^k, \ \text{Case 4: } \gamma_k = \gamma + 0.5^k.$$

Obviously, we have $\{|\gamma_k - \gamma|\}_{k\in\mathbb{N}} \in \ell_+^1$ for all the four cases. Problem (9.1) is considered again with $J$ the $\ell_1$, the $\ell_{1,2}$ and the nuclear norms. The comparison results are displayed in Figure 5. Table 2 shows the number of iteration $K$ needed for the identification of $\mathcal{M}_{x^\star}^J$.

For the stationary iteration, the local convergence rate of the 3 examples are,

$$\ell_1\text{-norm: } \rho = 0.9196, \quad \ell_{1,2}\text{-norm: } \rho = 0.9153, \quad \text{Nuclear norm: } \rho = 0.8904.$$

We can make the following observations from the comparison:

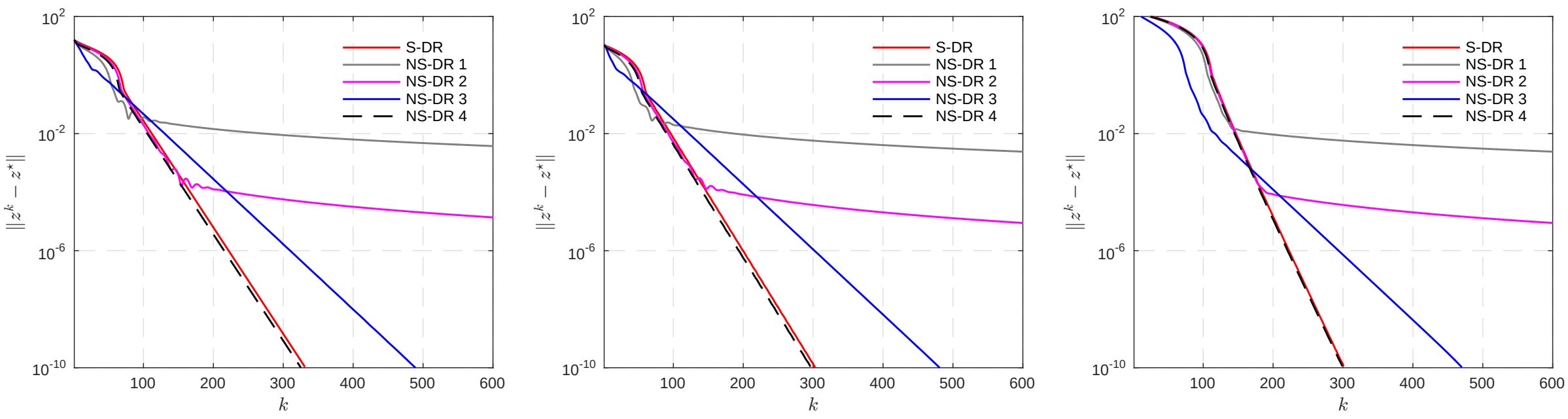

(a) $\ell_1$-norm (b) $\ell_{1,2}$-norm (c) Nuclear norm

Figure 5: Comparison between stationary ("S-DR") and non-stationary DR ("NS-DR X", X stands for Case X) when solving (9.1) with different functions $J$ in Table 1. (a) $\ell_1$-norm. (b) $\ell_{1,2}$-norm. (c) Nuclear norm.

Table 2: Number of iterations $K$ needed for the identification of $\mathcal{M}^J_{x^\star}$ for each tested case. "NS-DR X" stands for the non-stationary DR with choice of $\gamma_k$ as in Case X.

| | S-DR | NS-DR 1 | NS-DR 2 | NS-DR 3 | NS-DR 4 |
|---|---|---|---|---|---|
| $\ell_1$-norm | 81 | 64 | 77 | 244 | 75 |
| $\ell_{1,2}$-norm | 62 | 46 | 58 | 233 | 56 |
| Nuclear norm | 114 | 107 | 112 | 77 | 112 |

- The local convergence behaviour of the non-stationary iteration is no better than the stationary one which is in agreement with our analysis;
- As argued in Remark 5.7(ii), the convergence rate is eventually controlled by the error $|\gamma_k - \gamma|$, except for "Case 4", since $0.5$ is strictly smaller than the local linear rate of the stationary version (*i.e.* $|\gamma_k - \gamma| = o(\|z_k - z^\star\|)$);
- The non-stationary DR seems to generally lead to faster identification. But his is not a systematic behaviour as observed for instance for Case 3, where slower identification is obtained for the $\ell_1$ and the $\ell_{1,2}$ norms.

## Acknowledgments

This work has been partly supported by the European Research Council (ERC project SIGMA-Vision). JF was partly supported by Institut Universitaire de France. The authors would like to thank Russell Luke for helpful discussions.

## A Proof of Theorem 3.1

The following lemma is needed in the proof of Theorem 3.1.

**Lemma A.1.** *Suppose that conditions* (**H.2**) *and* (**H.3**) *hold, and that* $\gamma_k$ *is convergent. Then*

$$\lim_{k\to+\infty} \gamma_k = \gamma.$$

**Proof.** Since $\gamma_k$ is convergent, it has a unique cluster point, say $\lim_{k\to+\infty}\gamma_k = \gamma'$. It is then sufficient to show that $\gamma' = \gamma$. Suppose that $\gamma' \neq \gamma$. Fix some $\varepsilon \in ]0, |\gamma' - \gamma|[$. Thus, there exist an index $K > 0$ such that for all $k \geq K$,

$$|\gamma_k - \gamma'| < \varepsilon/2.$$

Therefore

$$|\gamma_k - \gamma| \geq |\gamma' - \gamma| - |\gamma_k - \gamma'| > \varepsilon/2.$$

It then follows that

$$\lambda_k(2-\lambda_k)\varepsilon \leq 2\lambda_k\varepsilon \leq 4\lambda_k|\gamma_k - \gamma|.$$

Denote $\overline{\tau} \stackrel{\text{def}}{=} \sup_{k\mathbb{N}} \lambda_k(2-\lambda_k)$ which is obviously positive and bounded since $\lambda_k \in [0,2]$. Summing both sides for $k \geq K$ we get

$$\begin{aligned}\varepsilon \textstyle\sum_{k\in\mathbb{N}} \lambda_k(2-\lambda_k) - K\overline{\tau} &\leq \varepsilon \textstyle\sum_{k=K}^{+\infty} \lambda_k(2-\lambda_k) \\ &\leq 4 \textstyle\sum_{k=K}^{\infty} \lambda_k|\gamma_k - \gamma| \\ &\leq 4 \textstyle\sum_{k\in\mathbb{N}} \lambda_k|\gamma_k - \gamma|,\end{aligned}$$

which, in view of (**H.3**), implies

$$\textstyle\sum_{k\in\mathbb{N}}\lambda_k(2-\lambda_k) \leq \varepsilon^{-1}(\lambda_k|\gamma_k-\gamma| + K\overline{\tau}) < +\infty,$$

which is a contradiction with (**H.2**). □

**Proof.** To prove our claim, we only need to check the conditions listed in [35, Theorem 4].

(i) As (**A.3**) assumes the set of minimizers of $(\mathcal{P})$ is nonempty, so is the set $\mathrm{Fix}(\mathscr{F}_\gamma)$, since the former is nothing but $\mathrm{prox}_{\gamma J}(\mathrm{Fix}(\mathscr{F}_\gamma))$ [4, Proposition 25.1(ii)].

(ii) Since $\mathscr{F}_{\gamma_k}$ is firmly nonexpansive by Lemma 2.4, $\mathscr{F}_{\gamma_k,\lambda_k}$ is $\frac{\lambda_k}{2}$-averaged nonexpansive, hence nonexpansive, owing to Lemma 2.3(iv).

(iii) Let $\rho \in [0,+\infty[$ and $z \in \mathbb{R}^n$ such that $\|z\| \leq \rho$, Then we have

$$\begin{aligned}(\mathscr{F}_{\gamma_k} - \mathscr{F}_\gamma)(z) &= \frac{\mathrm{rprox}_{\gamma_k G}\mathrm{o}\mathrm{rprox}_{\gamma_k J}}{2}(z) - \frac{\mathrm{rprox}_{\gamma G}\mathrm{o}\mathrm{rprox}_{\gamma J}}{2}(z) \\ &= \left(\frac{\mathrm{rprox}_{\gamma_k G}\mathrm{o}\mathrm{rprox}_{\gamma_k J}}{2}(z) - \frac{\mathrm{rprox}_{\gamma_k G}\mathrm{o}\mathrm{rprox}_{\gamma J}}{2}(z)\right) \\ &\quad - \left(\frac{\mathrm{rprox}_{\gamma G}\mathrm{o}\mathrm{rprox}_{\gamma J}}{2}(z) - \frac{\mathrm{rprox}_{\gamma_k G}\mathrm{o}\mathrm{rprox}_{\gamma J}}{2}(z)\right) \\ &= \left(\frac{\mathrm{rprox}_{\gamma_k G}\mathrm{o}\mathrm{rprox}_{\gamma_k J}}{2}(z) - \frac{\mathrm{rprox}_{\gamma_k G}\mathrm{o}\mathrm{rprox}_{\gamma J}}{2}(z)\right) \\ &\quad - \left(\mathrm{prox}_{\gamma G} \circ \mathrm{rprox}_{\gamma J}(z) - \mathrm{prox}_{\gamma_k G} \circ \mathrm{rprox}_{\gamma J}(z)\right).\end{aligned}$$

Thus, by virtue of Lemma 2.3(iii), we have

$$\begin{aligned}&\|(\mathscr{F}_{\gamma_k} - \mathscr{F}_\gamma)(z)\| \\ &\leq \|\mathrm{prox}_{\gamma_k J}(z) - \mathrm{prox}_{\gamma J}(z)\| + \|\mathrm{prox}_{\gamma_k G}(\mathrm{rprox}_{\gamma J}(z)) - \mathrm{prox}_{\gamma G}(\mathrm{rprox}_{\gamma J}(z))\|.\end{aligned}$$

Let's bound the first term. From the resolvent equation [9], and Lemma 2.3(i)(ii)(v), we have

$$\begin{aligned}\|\mathrm{prox}_{\gamma_k J}(z) - \mathrm{prox}_{\gamma J}(z)\| &= \|\mathrm{prox}_{\gamma_k J}(z) - \mathrm{prox}_{\gamma_k J}\big(\tfrac{\gamma_k}{\gamma}z + \big(1-\tfrac{\gamma_k}{\gamma}\big)\mathrm{prox}_{\gamma J}(z)\big)\| \\ &\leq \frac{|\gamma_k - \gamma|}{\gamma}\|(\mathrm{Id} - \mathrm{prox}_{\gamma J})(z)\| \leq \frac{|\gamma_k-\gamma|}{\gamma}(\rho + \|\mathrm{prox}_{\gamma J}(0)\|).\end{aligned} \tag{A.1}$$

With similar arguments, we also obtain

$$\|\mathrm{prox}_{\gamma_k G}(\mathrm{rprox}_{\gamma J}(z)) - \mathrm{prox}_{\gamma G}(\mathrm{rprox}_{\gamma J}(z))\| \leq \frac{|\gamma_k-\gamma|}{\gamma}(\rho + \|\mathrm{prox}_{\gamma G}(0)\| + 2\|\mathrm{prox}_{\gamma J}(0))\|). \tag{A.2}$$

Combining (A.1) and (A.2) leads to

$$\|(\mathscr{F}_{\gamma_k} - \mathscr{F}_\gamma)(z)\| \leq \frac{|\gamma_k - \gamma|}{\gamma}(2\rho + \|\mathrm{prox}_{\gamma G}(0)\| + 3\|\mathrm{prox}_{\gamma J}(0)\|), \tag{A.3}$$

whence we get

$$\begin{aligned}&\|(\mathscr{F}_{\gamma_k,\lambda_k} - \mathscr{F}_{\gamma,\lambda_k})(z)\| \\ &= \lambda_k\|(\mathscr{F}_{\gamma_k} - \mathscr{F}_\gamma)(z)\| \leq \lambda_k\frac{|\gamma_k - \gamma|}{\gamma}(2\rho + \|\mathrm{prox}_{\gamma G}(0)\| + 3\|\mathrm{prox}_{\gamma J}(0)\|).\end{aligned}$$

Therefore, from (**H.3**), we deduce that

$$\{\sup_{\|z\|\leq\rho} \|(\mathscr{F}_{\gamma_k,\lambda_k} - \mathscr{F}_{\gamma,\lambda_k})(z)\|\}_{k\in\mathbb{N}} \in \ell^1_+.$$

In other words, the non-stationary iteration (3.1) is a perturbed version of the stationary one (1.3) with an error term which is summable thanks to (**H.3**). The claim on the convergence of $z^\star$ follows by applying [13, Corollary 5.2]. Moreover, $x^\star \stackrel{\text{def}}{=} \mathrm{prox}_{\gamma J}(z^\star)$ is a solution of $(\mathcal{P})$. In turn, using nonexpansiveness of $\mathrm{prox}_{\gamma_k J}$ and (A.1), we have

$$\|x_k - x^\star\| \leq \|z_k - z^\star\| + \frac{|\gamma_k - \gamma|}{\gamma}(\|z^\star\| + \|\mathrm{prox}_{\gamma J}(0)\|),$$

and thus the right hand side goes to zero as $k \to +\infty$ as we are in finite dimension and since $\gamma_k \to \gamma$ owing to Lemma A.1. This entails that the shadow sequence $\{x_k\}_{k\in\mathbb{N}}$ also converges to $x^\star$. With similar arguments, we can also show that $\{v_k\}_{k\in\mathbb{N}}$ converges to $x^\star$ (using for instance (A.2) and nonexpansiveness of $\mathrm{prox}_{\gamma_k G}$). □

# B Proof of Section 4.1

## B.1 Proof of Theorem 4.3

**Proof.** By Theorem 3.1, we have the convergence of all the sequences generated by (1.5), that is

$$z_k \to z^\star \in \mathrm{Fix}(\mathscr{F}_{\gamma,\lambda}),\ x_k, v_k, \to x^\star = \mathrm{prox}_{\gamma J}(z^\star) \in \mathrm{Argmin}(G + J).$$

The nondegeneracy condition (ND) is equivalent to

$$\frac{x^\star - z^\star}{\gamma} \in \mathrm{ri}\big(\partial G(x^\star)\big) \ \text{ and } \ \frac{z^\star - x^\star}{\gamma} \in \mathrm{ri}\big(\partial J(x^\star)\big). \tag{B.1}$$

(i) The update of $x_{k+1}$ and $v_{k+1}$ in iteration (1.5) is equivalent to the monotone inclusions

$$\frac{2x_k - z_k - v_{k+1}}{\gamma_k} \in \partial G(v_{k+1}) \ \text{ and } \ \frac{z_k - x_k}{\gamma_k} \in \partial J(x_k).$$

It then follows that

$$\begin{aligned}&\mathrm{dist}\big(\frac{x^\star - z^\star}{\gamma}, \partial G(v_{k+1})\big) \\ &\leq \|\frac{x^\star - z^\star}{\gamma} - \frac{2x_k - z_k - v_{k+1}}{\gamma_k}\| \\ &= \|\frac{(\gamma_k - \gamma)(x^\star - z^\star)}{\gamma\gamma_k} + \frac{x^\star - z^\star}{\gamma_k} - \frac{2x_k - z_k - v_{k+1}}{\gamma_k}\| \\ &\leq \frac{|\gamma_k - \gamma|}{\gamma\underline{\gamma}}\|(\mathrm{Id} - \mathrm{prox}_{\gamma J})(z^\star)\| + \frac{1}{\underline{\gamma}}\|(z_k - z^\star) - 2(x_k - x^\star) + (v_{k+1} - x^\star)\| \\ &\leq \frac{|\gamma_k - \gamma|}{\gamma\underline{\gamma}}(\|z^\star\| + \mathrm{prox}_{\gamma J}(0)) + \frac{1}{\underline{\gamma}}(\|z_k - z^\star\| + 2\|x_k - x^\star\| + \|v_{k+1} - x^\star\|),\end{aligned}$$

and the right hand side converges to 0 in view of Theorem 3.1 and Lemma A.1. Similarly, we have

$$\begin{aligned}
\mathrm{dist}\big(\frac{z^\star - x^\star}{\gamma}, \partial J(x_k)\big) &\leq \|\frac{z^\star - x^\star}{\gamma} - \frac{z_k - x_k}{\gamma_k}\| \\
&= \|\frac{(\gamma_k - \gamma)(z^\star - x^\star)}{\gamma\gamma_k} + \frac{z^\star - x^\star}{\gamma_k} - \frac{z_k - x_k}{\gamma_k}\| \\
&\leq \frac{|\gamma_k - \gamma|}{\gamma\underline{\gamma}}(\|z^\star\| + \mathrm{prox}_{\gamma J}(0)) + \frac{1}{\underline{\gamma}}(\|z_k - z^\star\| + \|x_k - x^\star\|) \to 0.
\end{aligned}$$

By assumption, $G, J \in \Gamma_0(\mathbb{R}^n)$, hence are subdifferentially continuous at every point in their respective domains [44, Example 13.30], and in particular at $x^\star$. It then follows that $G(v_k) \to G(x^\star)$ and $J(x_k) \to J(x^\star)$. Altogether, this shows that the conditions of [24, Theorem 5.3] are fulfilled for $G$ and $J$, and the finite identification claim follows.

(ii) (a) In this case, $\mathcal{M}^J_{x^\star}$ is an affine subspace, *i.e.* $\mathcal{M}^J_{x^\star} = x^\star + T^J_{x^\star}$. Since $J$ is partly smooth at $x^\star$ relative to $\mathcal{M}^J_{x^\star}$, the sharpness property holds at all nearby points in $\mathcal{M}^J_{x^\star}$ [31, Proposition 2.10]. Thus for $k$ large enough, *i.e.* $x_k$ sufficiently close to $x^\star$ on $\mathcal{M}^J_{x^\star}$, we have indeed $\mathcal{T}_{x_k}(\mathcal{M}^J_{x^\star}) = T^J_{x^\star} = T^J_{x_k}$ as claimed.

(b) Similar to (ii)(a).

(c) It is immediate to verify that a locally polyhedral function around $x^\star$ is indeed partly smooth relative to the affine subspace $x^\star + T^J_{x^\star}$, and thus, the first claim follows from (ii)(a). For the rest, it is sufficient to observe that by polyhedrality, for any $x \in \mathcal{M}^J_{x^\star}$ near $x^\star$, $\partial J(x) = \partial J(x^\star)$. Therefore, combining local normal sharpness [31, Proposition 2.10] and Lemma 4.2 yields the second conclusion.

(d) Similar to (ii)(c). □

## B.2 Proof of Proposition 4.5

**Proof.** From (3.1), we have

$$z_{k+1} = \mathscr{F}_{\gamma,\lambda_k}(z_k) + e_k$$

where $\{\|e_k\|\}_{k\in\mathbb{N}} = \{O(\lambda_k|\gamma_k - \gamma|)\}_{k\in\mathbb{N}} \in \ell^1_+$ (see the proof of Theorem 3.1). Since $\mathscr{F}_{\gamma_k}$ is firmly non-expansive by Lemma 2.4, $\mathscr{F}_{\gamma,\lambda_k}$ is $\frac{\lambda_k}{2}$-averaged non-expansive owing to Lemma 2.3(iv). Thus arguing as in the proof of [13, Theorem 3.1], we have

$$\begin{aligned}
\|z_k - z^\star\|^2 &\leq \|\mathscr{F}_{\gamma,\lambda_k}(z_{k-1}) - \mathscr{F}_{\gamma,\lambda_k}(z^\star)\|^2 + C\|e_{k-1}\| \\
&\leq \|\mathscr{F}_{\gamma,\lambda_k}(z_{k-1}) - \mathscr{F}_{\gamma,\lambda_k}(z^\star)\|^2 - \tfrac{2-\lambda_{k-1}}{\lambda_{k-1}}\|z_k - z_{k-1}\|^2 + C\|e_{k-1}\| \\
&\leq \|z_{k-1} - z^\star\|^2 - \tau_{k-1}\|v_k - x_{k-1}\|^2 + C\|e_{k-1}\|,
\end{aligned}$$

where $C < +\infty$ by boundedness of $z_k$ and $e_k$. Let $g_k = (z_{k-1} - x_{k-1})/\gamma_{k-1}$ and $h_k = (2x_{k-1} - z_{k-1} - v_k)/\gamma_{k-1}$. By definition, we have $(g_k, h_k) \in \partial J(x_{k-1}) \times \partial G(v_k)$. Suppose that neither $\mathcal{M}^J_{x^\star}$ nor $\mathcal{M}^G_{x^\star}$ have been identified at iteration $k$. That is $x_{k-1} \notin \mathcal{M}^J_{x^\star}$ and $v_k \notin \mathcal{M}^G_{x^\star}$, and by assumption, $g_k \in \mathrm{rbd}(\partial J(x^\star))$ and $h_k \in \mathrm{rbd}(\partial G(x^\star))$, which implies that $g_k + h_k = (v_k - x_{k-1})/\gamma_{k-1} \in \mathrm{rbd}(\partial J(x^\star)) + \mathrm{rbd}(\partial G(x^\star))$. Thus, the above inequality becomes

$$\begin{aligned}
\|z_k - z^\star\|^2 &\leq \|z_{k-1} - z^\star\|^2 - \gamma_{k-1}^2\tau_{k-1}\mathrm{dist}\big(0, \mathrm{rbd}(\partial J(x^\star)) + \mathrm{rbd}(\partial G(x^\star))\big)^2 + C\|e_{k-1}\| \\
&\leq \|z_{k-1} - z^\star\|^2 - \gamma_{k-1}^2\tau_{k-1}\mathrm{dist}\big(0, \mathrm{rbd}(\partial J(x^\star) + \partial G(x^\star))\big)^2 + C\|e_{k-1}\| \\
&\leq \|z_0 - z^\star\|^2 - k\underline{\gamma}^2\underline{\tau}\mathrm{dist}\big(0, \mathrm{rbd}(\partial J(x^\star) + \partial G(x^\star))\big)^2 + O\big(\textstyle\sum_{k\in\mathbb{N}}\lambda_k|\gamma_k - \gamma|\big),
\end{aligned}$$

and $\mathrm{dist}\big(0, \mathrm{rbd}(\partial J(x^\star) + \partial G(x^\star))\big) > 0$ owing to condition (ND). Taking $k$ as the largest integer such that the bound in the right hand is positive, we deduce that the number of iterations where both $\mathcal{M}^J_{x^\star}$ and $\mathcal{M}^G_{x^\star}$ have not been identified yet does not exceed the claimed bound (4.1). Thus finite identification necessarily occurs at some $k$ larger than this bound. □

## C Proofs of Section 5

### C.1 Riemannian Geometry

Let $\mathcal{M}$ be a $C^2$-smooth embedded submanifold of $\mathbb{R}^n$ around a point $x$. With some abuse of terminology, we shall state $C^2$-manifold instead of $C^2$-smooth embedded submanifold of $\mathbb{R}^n$. The natural embedding of a submanifold $\mathcal{M}$ into $\mathbb{R}^n$ permits to define a Riemannian structure and to introduce geodesics on $\mathcal{M}$, and we simply say $\mathcal{M}$ is a Riemannian manifold. We denote respectively $\mathcal{T}_{\mathcal{M}}(x)$ and $\mathcal{N}_{\mathcal{M}}(x)$ the tangent and normal space of $\mathcal{M}$ at point near $x$ in $\mathcal{M}$.

**Exponential map** Geodesics generalize the concept of straight lines in $\mathbb{R}^n$, preserving the zero acceleration characteristic, to manifolds. Roughly speaking, a geodesic is locally the shortest path between two points on $\mathcal{M}$. We denote by $\mathfrak{g}(t; x, h)$ the value at $t \in \mathbb{R}$ of the geodesic starting at $\mathfrak{g}(0; x, h) = x \in \mathcal{M}$ with velocity $\dot{\mathfrak{g}}(t; x, h) = \frac{d\mathfrak{g}}{dt}(t; x, h) = h \in \mathcal{T}_{\mathcal{M}}(x)$ (which is uniquely defined). For every $h \in \mathcal{T}_{\mathcal{M}}(x)$, there exists an interval $I$ around 0 and a unique geodesic $\mathfrak{g}(t; x, h) : I \to \mathcal{M}$ such that $\mathfrak{g}(0; x, h) = x$ and $\dot{\mathfrak{g}}(0; x, h) = h$. The mapping

$$\mathrm{Exp}_x : \mathcal{T}_{\mathcal{M}}(x) \to \mathcal{M}, \;\; h \mapsto \mathrm{Exp}_x(h) = \mathfrak{g}(1; x, h),$$

is called *Exponential map*. Given $x, x' \in \mathcal{M}$, the direction $h \in \mathcal{T}_{\mathcal{M}}(x)$ we are interested in is such that

$$\mathrm{Exp}_x(h) = x' = \mathfrak{g}(1; x, h).$$

**Parallel translation** Given two points $x, x' \in \mathcal{M}$, let $\mathcal{T}_{\mathcal{M}}(x), \mathcal{T}_{\mathcal{M}}(x')$ be their corresponding tangent spaces. Define

$$\tau : \mathcal{T}_{\mathcal{M}}(x) \to \mathcal{T}_{\mathcal{M}}(x'),$$

the parallel translation along the unique geodesic joining $x$ to $x'$, which is isomorphism and isometry w.r.t. the Riemannian metric.

**Riemannian gradient and Hessian** For a vector $v \in \mathcal{N}_{\mathcal{M}}(x)$, the Weingarten map of $\mathcal{M}$ at $x$ is the operator $\mathfrak{W}_x(\cdot, v) : \mathcal{T}_{\mathcal{M}}(x) \to \mathcal{T}_{\mathcal{M}}(x)$ defined by

$$\mathfrak{W}_x(\cdot, v) = -\mathrm{P}_{\mathcal{T}_{\mathcal{M}}(x)} \mathrm{d}V[h],$$

where $V$ is any local extension of $v$ to a normal vector field on $\mathcal{M}$. The definition is independent of the choice of the extension $V$, and $\mathfrak{W}_x(\cdot, v)$ is a symmetric linear operator which is closely tied to the second fundamental form of $\mathcal{M}$, see [10, Proposition II.2.1].

Let $G$ be a real-valued function which is $C^2$ along the $\mathcal{M}$ around $x$. The covariant gradient of $G$ at $x' \in \mathcal{M}$ is the vector $\nabla_{\mathcal{M}} G(x') \in \mathcal{T}_{\mathcal{M}}(x')$ defined by

$$\langle \nabla_{\mathcal{M}} G(x'), \, h \rangle = \frac{d}{dt} G\big(\mathrm{P}_{\mathcal{M}}(x' + th)\big)\big|_{t=0}, \;\; \forall h \in \mathcal{T}_{\mathcal{M}}(x'),$$

where $\mathrm{P}_{\mathcal{M}}$ is the projection operator onto $\mathcal{M}$. The covariant Hessian of $G$ at $x'$ is the symmetric linear mapping $\nabla^2_{\mathcal{M}} G(x')$ from $\mathcal{T}_{\mathcal{M}}(x')$ to itself which is defined as

$$\langle \nabla^2_{\mathcal{M}} G(x') h, \, h \rangle = \frac{d^2}{dt^2} G\big(\mathrm{P}_{\mathcal{M}}(x' + th)\big)\big|_{t=0}, \;\; \forall h \in \mathcal{T}_{\mathcal{M}}(x'). \tag{C.1}$$

This definition agrees with the usual definition using geodesics or connections [39]. Now assume that $\mathcal{M}$ is a Riemannian embedded submanifold of $\mathbb{R}^n$, and that a function $G$ has a $C^2$-smooth restriction on $\mathcal{M}$. This can be characterized by the existence of a $C^2$-smooth extension (representative) of $G$, *i.e.* a $C^2$-smooth function $\widetilde{G}$ on $\mathbb{R}^n$ such that $\widetilde{G}$ agrees with $G$ on $\mathcal{M}$. Thus, the Riemannian gradient $\nabla_{\mathcal{M}} G(x')$ is also given by

$$\nabla_{\mathcal{M}} G(x') = \mathrm{P}_{\mathcal{T}_{\mathcal{M}}(x')} \nabla \widetilde{G}(x'), \tag{C.2}$$

and $\forall h \in \mathcal{T}_{\mathcal{M}}(x')$, the Riemannian Hessian reads

$$\begin{aligned}\nabla^2_{\mathcal{M}} G(x') h &= \mathrm{P}_{\mathcal{T}_{\mathcal{M}}(x')} \mathrm{d}(\nabla_{\mathcal{M}} G)(x')[h] = \mathrm{P}_{\mathcal{T}_{\mathcal{M}}(x')} \mathrm{d}\big(x' \mapsto \mathrm{P}_{\mathcal{T}_{\mathcal{M}}(x')} \nabla_{\mathcal{M}} \widetilde{G}\big)[h] \\ &= \mathrm{P}_{\mathcal{T}_{\mathcal{M}}(x')} \nabla^2 \widetilde{G}(x') h + \mathfrak{W}_{x'}\big(h, \mathrm{P}_{\mathcal{N}_{\mathcal{M}}(x')} \nabla \widetilde{G}(x')\big),\end{aligned} \tag{C.3}$$

where the last equality comes from [1, Theorem 1]. When $\mathcal{M}$ is an affine or linear subspace of $\mathbb{R}^n$, then obviously $\mathcal{M} = x + \mathcal{T}_{\mathcal{M}}(x)$, and $\mathfrak{W}_{x'}(h, \mathrm{P}_{\mathcal{N}_{\mathcal{M}}(x')} \nabla \widetilde{G}(x')) = 0$, hence (C.3) reduces to

$$\nabla^2_{\mathcal{M}} G(x') = \mathrm{P}_{\mathcal{T}_{\mathcal{M}}(x')} \nabla^2 \widetilde{G}(x') \mathrm{P}_{\mathcal{T}_{\mathcal{M}}(x')}.$$

See [29, 10] for more materials on differential and Riemannian manifolds.

We have the following proposition characterising the parallel translation and the Riemannian Hessian of two close points in $\mathcal{M}$.

**Lemma C.1.** *Let $x, x'$ be two close points in $\mathcal{M}$, denote $\mathcal{T}_{\mathcal{M}}(x), \mathcal{T}_{\mathcal{M}}(x')$ be the tangent spaces of $\mathcal{M}$ at $x, x'$ respectively, and $\tau : \mathcal{T}_{\mathcal{M}}(x') \to \mathcal{T}_{\mathcal{M}}(x)$ be the parallel translation along the unique geodesic joining from $x$ to $x'$, then for the parallel translation we have, given any bounded vector $v \in \mathbb{R}^n$*

$$(\tau \mathrm{P}_{\mathcal{T}_{\mathcal{M}}(x')} - \mathrm{P}_{\mathcal{T}_{\mathcal{M}}(x)}) v = o(v). \tag{C.4}$$

*The Riemannian Taylor expansion of $J \in C^2(\mathcal{M})$ at $x$ for $x'$ reads,*

$$\tau \nabla_{\mathcal{M}} J(x') = \nabla_{\mathcal{M}} J(x) + \nabla^2_{\mathcal{M}} J(x) \mathrm{P}_{\mathcal{T}_{\mathcal{M}}(x)}(x' - x) + o(x' - x). \tag{C.5}$$

**Proof.** See [34, Lemma B.1 and B.2]. □

### C.2 Proof of Proposition 5.2

**Proof.** Since $W_{\overline{G}}, W_{\overline{J}}$ are both firmly non-expansive by Lemma 5.1, it follows from [4, Example 4.7] that $M_{\overline{G}}$ and $M_{\overline{J}}$ are firmly non-expansive. As a result, $M$ is firmly non-expansive [4, Proposition 4.21(i)-(ii)], and equivalently that $M_\lambda$ is $\frac{\lambda}{2}$-averaged by Lemma 2.3(i)⇔(iv).

Under the assumptions of Theorem 4.3, there exists $K \in \mathbb{N}$ large enough such that for all $k \geq K$, $(x_k, v_k) \in \mathcal{M}^J_{x^\star} \times \mathcal{M}^G_{x^\star}$. Denote $T^J_{x_k}$ and $T^J_{x^\star}$ be the tangent spaces corresponding to $x_k$ and $x^\star \in \mathcal{M}^J_{x^\star}$, and similarly $T^G_{x_k}$ and $T^G_{x^\star}$ the tangent spaces corresponding to $v_k$ and $x^\star \in \mathcal{M}^G_{x^\star}$. Denote $\tau^J_k : T^J_{x_k} \to T^J_{x^\star}$ (resp. $\tau^G_k : T^G_{v_k} \to T^G_{x^\star}$) the parallel translation along the unique geodesic on $\mathcal{M}^J_{x^\star}$ (resp. $\mathcal{M}^G_{x^\star}$) joining $x_k$ to $x^\star$ (resp. $v_k$ to $x^\star$).

From (1.5), for $x_k$, we have

$$\begin{cases} x_k = \mathrm{prox}_{\gamma_k J}(z_k), \\ x^\star = \mathrm{prox}_{\gamma J}(z^\star), \end{cases} \Longleftrightarrow \begin{cases} z_k - x_k \in \gamma_k \partial J(x_k), \\ z^\star - x^\star \in \gamma \partial J(x^\star). \end{cases}$$

Projecting on the corresponding tangent spaces, using Lemma 4.2, and applying the parallel translation operator $\tau^J_k$ leads to

$$\begin{aligned}\gamma_k \tau^J_k \nabla_{\mathcal{M}^J_{x^\star}} J(x_k) &= \tau^J_k \mathrm{P}_{T^J_{x_k}}(z_k - x_k) = \mathrm{P}_{T^J_{x^\star}}(z_k - x_k) + \big(\tau^J_k \mathrm{P}_{T^J_{x_k}} - \mathrm{P}_{T^J_{x^\star}}\big)(z_k - x_k), \\ \gamma \nabla_{\mathcal{M}^J_{x^\star}} J(x^\star) &= \mathrm{P}_{T^J_{x^\star}}(z^\star - x^\star).\end{aligned}$$

We then obtain

$$
\begin{aligned}
&\gamma_k \tau_k^J \nabla_{\mathcal{M}_{x^\star}^J} J(x_k) - \gamma \nabla_{\mathcal{M}_{x^\star}^J} J(x^\star) \\
&= \gamma \tau_k^J \nabla_{\mathcal{M}_{x^\star}^J} J(x_k) - \gamma \nabla_{\mathcal{M}_{x^\star}^J} J(x^\star) + (\gamma_k - \gamma)\tau_k^J \nabla_{\mathcal{M}_{x^\star}^J} J(x_k) \\
&= \mathrm{P}_{T_{x^\star}^J}\big((z_k - z^\star) - (x_k - x^\star)\big) \\
&\quad + \underbrace{\big(\tau_k^J \mathrm{P}_{T_{x_k}^J} - \mathrm{P}_{T_{x^\star}^J}\big)(z_k - x_k - z^\star + x^\star)}_{\textbf{Term 1}} + \underbrace{\big(\tau_k^J \mathrm{P}_{T_{x_k}^J} - \mathrm{P}_{T_{x^\star}^J}\big)(z^\star - x^\star)}_{\textbf{Term 2}}.
\end{aligned}
\tag{C.6}
$$

For $(\gamma_k - \gamma)\tau_k^J \nabla_{\mathcal{M}_{x^\star}^J} J(x_k)$, since the Riemannian gradient $\nabla_{\mathcal{M}_{x^\star}^J} J(x_k)$ is single-valued and bounded on bounded sets, we have

$$
\|(\gamma_k - \gamma)\tau_k^J \nabla_{\mathcal{M}_{x^\star}^J} J(x_k)\| = O(|\gamma_k - \gamma|). \tag{C.7}
$$

Combining (A.1) and (C.4), we have for **Term 1**

$$
(\tau_k^J \mathrm{P}_{T_{x_k}^J} - \mathrm{P}_{T_{x^\star}^J})(z_k - x_k - z^\star + x^\star) = o(\|z_k - z^\star\|) + o(|\gamma_k - \gamma|). \tag{C.8}
$$

As far as **Term 2** is concerned, with (5.1), (A.1) and the Riemannian Taylor expansion (C.5), we have

$$
\begin{aligned}
&\gamma \tau_k^J \nabla_{\mathcal{M}_{x^\star}^J} J(x_k) - \gamma \nabla_{\mathcal{M}_{x^\star}^J} J(x^\star) - \big(\tau_k^J \mathrm{P}_{T_{x_k}^J} - \mathrm{P}_{T_{x^\star}^J}\big)(z^\star - x^\star) \\
&= \tau_k^J\big(\gamma \nabla_{\mathcal{M}_{x^\star}^J} J(x_k) - \mathrm{P}_{T_{x_k}^J}(z^\star - x^\star)\big) - \big(\gamma \nabla_{\mathcal{M}_{x^\star}^J} J(x^\star) - \mathrm{P}_{T_{x^\star}^J}(z^\star - x^\star)\big) \\
&= \tau_k^J \nabla_{\mathcal{M}_{x^\star}^J} \bar{J}(x_k) - \nabla_{\mathcal{M}_{x^\star}^J} \bar{J}(x^\star) = \mathrm{P}_{T_{x^\star}^J} \nabla^2_{\mathcal{M}_{x^\star}^J} \bar{J}(x^\star) \mathrm{P}_{T_{x^\star}^J}(x_k - x^\star) + o(x_k - x^\star) \\
&= \mathrm{P}_{T_{x^\star}^J} \nabla^2_{\mathcal{M}_{x^\star}^J} \bar{J}(x^\star) \mathrm{P}_{T_{x^\star}^J}(x_k - x^\star) + o(\|z_k - z^\star\|) + o(|\gamma_k - \gamma|).
\end{aligned}
\tag{C.9}
$$

Therefore, inserting (C.7), (C.8) and (C.9) into (C.6), we obtain

$$
\begin{aligned}
& H_{\bar{J}}(x_k - x^\star) = \mathrm{P}_{T_{x^\star}^J}(z_k - z^\star) - \mathrm{P}_{T_{x^\star}^J}(x_k - x^\star) + o(\|z_k - z^\star\|) + O(|\gamma_k - \gamma|) \\
\Rightarrow\ & (\mathrm{Id} + H_{\bar{J}})\mathrm{P}_{T_{x^\star}^J}(x_k - x^\star) = \mathrm{P}_{T_{x^\star}^J}(z_k - z^\star) + o(\|z_k - z^\star\|) + O(|\gamma_k - \gamma|) \\
\Rightarrow\ & \mathrm{P}_{T_{x^\star}^J}(x_k - x^\star) = W_{\bar{J}} \mathrm{P}_{T_{x^\star}^J}(z_k - z^\star) + o(\|z_k - z^\star\|) + O(|\gamma_k - \gamma|) \\
\Rightarrow\ & \mathrm{P}_{T_{x^\star}^J}(x_k - x^\star) = \mathrm{P}_{T_{x^\star}^J} W_{\bar{J}} \mathrm{P}_{T_{x^\star}^J}(z_k - z^\star) + o(\|z_k - z^\star\|) + O(|\gamma_k - \gamma|) \\
\Rightarrow\ & x_k - x^\star = M_{\bar{J}}(z_k - z^\star) + o(\|z_k - z^\star\|) + O(|\gamma_k - \gamma|),
\end{aligned}
\tag{C.10}
$$

where we used the fact that $x_k - x^\star = \mathrm{P}_{T_{x^\star}^J}(x_k - x^\star) + o(x_k - x^\star)$ [33, Lemma 5.1].

Similarly for $v_{k+1}$, we have

$$
\begin{cases} v_{k+1} = \mathrm{prox}_{\gamma_k G}(2x_k - z_k), \\ x^\star = \mathrm{prox}_{\gamma G}(2x^\star - z^\star), \end{cases} \iff \begin{cases} 2x_k - z_k - v_{k+1} \in \gamma \partial J(v_{k+1}), \\ 2x^\star - z^\star - x^\star \in \gamma \partial J(x^\star). \end{cases}
$$

Upon projecting onto the corresponding tangent spaces and applying the parallel translation $\tau_{k+1}^G$, we get

$$
\begin{aligned}
\gamma_k \tau_{k+1}^G \nabla_{\mathcal{M}_{x^\star}^G} G(v_{k+1}) &= \tau_{k+1}^G \mathrm{P}_{T_{v_{k+1}}^G}(2x_k - z_k - v_{k+1}) \\
&= \mathrm{P}_{T_{x^\star}^G}(2x_k - z_k - v_{k+1}) + \big(\tau_{k+1}^G \mathrm{P}_{T_{v_{k+1}}^G} - \mathrm{P}_{T_{x^\star}^G}\big)(2x_k - z_k - v_{k+1}), \\
\gamma \nabla_{\mathcal{M}_{x^\star}^G} G(x^\star) &= \mathrm{P}_{T_{x^\star}^G}(2x^\star - z^\star - x^\star).
\end{aligned}
$$

Substracting both equations, we obtain

$$\begin{aligned}
&\gamma\tau^G_{k+1}\nabla_{\mathcal{M}^G_{x^\star}} G(v_{k+1}) - \gamma\nabla_{\mathcal{M}^G_{x^\star}} G(x^\star) \\
&= \gamma\tau^G_{k+1}\nabla_{\mathcal{M}^G_{x^\star}} G(v_{k+1}) - \gamma\nabla_{\mathcal{M}^G_{x^\star}} G(x^\star) + (\gamma_k - \gamma)\tau^G_{k+1}\nabla_{\mathcal{M}^G_{x^\star}} G(v_{k+1}) \\
&= \mathrm{P}_{T^G_{x^\star}}\big((2x_k - z_k - v_{k+1}) - (2x^\star - z^\star - x^\star)\big) + \underbrace{\big(\tau^G_{k+1}\mathrm{P}_{T^G_{v_{k+1}}} - \mathrm{P}_{T^G_{x^\star}}\big)(x^\star - z^\star)}_{\textbf{Term 4}} \\
&\quad + \underbrace{\big(\tau^G_{k+1}\mathrm{P}_{T^G_{v_{k+1}}} - \mathrm{P}_{T^G_{x^\star}}\big)\big((2x_k - z_k - v_{k+1}) - (2x^\star - z^\star - x^\star)\big)}_{\textbf{Term 3}}.
\end{aligned} \tag{C.11}$$

As for (C.7), we have

$$\|(\gamma_k - \gamma)\tau^G_{k+1}\nabla_{\mathcal{M}^G_{x^\star}} G(v_{k+1})\| = O(|\gamma_k - \gamma|). \tag{C.12}$$

With similar arguments to those used for **Term 1**, we have **Term 3** $= o(\|z_k - z^\star\|) + o(|\gamma_k - \gamma|)$. Moreover, similarly to (C.9), we have for **Term 4**,

$$\begin{aligned}
&\gamma\tau^G_{k+1}\nabla_{\mathcal{M}^G_{x^\star}} G(v_{k+1}) - \gamma\nabla_{\mathcal{M}^G_{x^\star}} G(x^\star) - \big(\tau^G_{k+1}\mathrm{P}_{T^G_{v_{k+1}}} - \mathrm{P}_{T^G_{x^\star}}\big)(x^\star - z^\star) \\
&= \mathrm{P}_{T^G_{x^\star}}\nabla^2_{\mathcal{M}^G_{x^\star}}\overline{G}(x^\star)\mathrm{P}_{T^G_{x^\star}}(v_{k+1} - x^\star) + o(z_k - x^\star) + o(|\gamma_k - \gamma|).
\end{aligned} \tag{C.13}$$

Then for (C.11) we have,

$$\begin{aligned}
&& H_{\overline{G}}(v_{k+1} - x^\star) &= 2\mathrm{P}_{T^G_{x^\star}}(x_k - x^\star) - \mathrm{P}_{T^G_{x^\star}}(z_k - z^\star) - \mathrm{P}_{T^G_{x^\star}}(v_{k+1} - x^\star) \\
&& &\quad + o(\|z_k - z^\star\|) + O(|\gamma_k - \gamma|) \\
\Rightarrow && (\mathrm{Id} + H_{\overline{G}})\mathrm{P}_{T^G_{x^\star}}(v_{k+1} - x^\star) &= 2\mathrm{P}_{T^G_{x^\star}}(x_k - x^\star) - \mathrm{P}_{T^G_{x^\star}}(z_k - z^\star) + o(\|z_k - z^\star\|) + O(|\gamma_k - \gamma|) \\
\Rightarrow && \mathrm{P}_{T^G_{x^\star}}(v_{k+1} - x^\star) &= 2M_{\overline{G}}M_{\overline{J}}(z_k - z^\star) - M_{\overline{G}}(z_k - z^\star) + o(\|z_k - z^\star\|) + O(|\gamma_k - \gamma|) \\
\Rightarrow && v_{k+1} - x^\star &= 2M_{\overline{G}}M_{\overline{J}}(z_k - z^\star) - M_{\overline{G}}(z_k - z^\star) + o(\|z_k - z^\star\|) + O(|\gamma_k - \gamma|),
\end{aligned} \tag{C.14}$$

where $v_{k+1} - x^\star = \mathrm{P}_{T^G_{x^\star}}(v_{k+1} - x^\star) + o(v_{k+1} - x^\star)$ is applied again [33, Lemma 5.1].

Summing up (C.10) and (C.14), we get

$$\begin{aligned}
&(z_k + v_{k+1} - x_k) - (z^\star + x^\star - x^\star) \\
&= (z_k - z^\star) + (v_{k+1} - x^\star) - (x_k - x^\star) \\
&= (\mathrm{Id} + 2M_{\overline{G}}M_{\overline{J}} - M_{\overline{G}} - M_{\overline{J}})(z_k - z^\star) + o(\|z_k - z^\star\|) + O(|\gamma_k - \gamma|) \\
&= M(z_k - z^\star) + o(\|z_k - z^\star\|) + O(|\gamma_k - \gamma|).
\end{aligned}$$

Hence for the relaxed DR iteration, we have

$$\begin{aligned}
z_{k+1} - z^\star &= (1 - \lambda_k)(z_k - z^\star) + \lambda_k\big((z_k + v_{k+1} - x_k) - (z^\star + x^\star - x^\star)\big) \\
&= (1 - \lambda_k)(z_k - z^\star) + \lambda_k M(z_k - z^\star) + o(\|z_k - z^\star\|) + \phi_k \\
&= M_\lambda(z_k - z^\star) - (\lambda_k - \lambda)(\mathrm{Id} - M)(z_k - z^\star) + o(\|z_k - z^\star\|) + \phi_k
\end{aligned}$$

Since $\mathrm{Id} - M$ is also (firmly) non-expansive (Lemma 2.3(ii)) and $\lambda_k \to \lambda \in ]0, 2[$, we thus get

$$\lim_{k\to\infty}\frac{\|(\lambda_k - \lambda)(\mathrm{Id} - M)(z_k - z^\star)\|}{\|z_k - z^\star\|} = \lim_{k\to\infty}\frac{|\lambda_k - \lambda|\|(\mathrm{Id} - M)(z_k - z^\star)\|}{\|z_k - z^\star\|} \le \lim_{k\to\infty}|\lambda_k - \lambda| = 0,$$

which means that

$$z_{k+1} - z^\star = M_\lambda(z_k - z^\star) + \psi_k + \phi_k,$$

and the claimed result is obtained. □

## C.3 Proof of Lemma 5.4

**Proof.**

(i) Since $M$ is firmly non-expansive and $M_\lambda$ is $\frac{\lambda}{2}$-averaged by Proposition 5.2, we deduce from [4, Proposition 5.15] that $M$ and $M_\lambda$ are convergent, and their limit is $M_\lambda^\infty = \mathrm{P}_{\mathrm{Fix}(M_\lambda)} = \mathrm{P}_{\mathrm{Fix}(M)} = M^\infty$ [3, Corollary 2.7(ii)]. Moreover, $M_\lambda^k - M^\infty = (M_\lambda - M^\infty)^k, \forall k \in \mathbb{N}$, and $\rho(M_\lambda - M^\infty) < 1$ by [3, Theorem 2.12]. It is also immediate to see that

$$\mathrm{Fix}(M) = \ker\big(M_{\overline{G}}(\mathrm{Id} - M_{\overline{J}}) + (\mathrm{Id} - M_{\overline{G}})M_{\overline{J}}\big).$$

Observe that

$$\begin{aligned} \mathrm{span}(M_{\overline{J}}) \subseteq T_{x^\star}^J \ \text{ and } \ \mathrm{span}(M_{\overline{G}}) \subseteq T_{x^\star}^G, \\ \ker\big(\mathrm{Id} - M_{\overline{G}}\big) \subseteq T_{x^\star}^G \ \text{ and } \ \ker(M_{\overline{G}}) = S_{x^\star}^G, \\ \mathrm{span}\big((\mathrm{Id} - M_{\overline{G}})M_{\overline{J}}\big) \subseteq \mathrm{span}(\mathrm{Id} - M_{\overline{G}}) \ \text{ and } \ \mathrm{span}\big(M_{\overline{G}}(\mathrm{Id} - M_{\overline{J}})\big) \subseteq T_{x^\star}^G, \end{aligned}$$

where we used the fact that $W_{\overline{G}}$ and $W_{\overline{J}}$ are positive definite. Therefore, $M_\lambda^\infty = 0$, if and only if, $\mathrm{Fix}(M) = \{0\}$, and for this to hold true, it is sufficient that

$$\begin{aligned} \mathrm{span}(M_{\overline{J}}) \cap \ker(\mathrm{Id} - M_{\overline{G}}) \subseteq T_{x^\star}^J \cap T_{x^\star}^G = \{0\}, \\ \mathrm{span}(\mathrm{Id} - M_{\overline{J}}) \cap \ker(M_{\overline{G}}) = \mathrm{span}(\mathrm{Id} - M_{\overline{J}}) \cap S_{x^\star}^G = \{0\}, \\ \mathrm{span}\big((\mathrm{Id} - M_{\overline{G}})M_{\overline{J}}\big) \cap \mathrm{span}\big(M_{\overline{G}}(\mathrm{Id} - M_{\overline{J}})\big) \subseteq \mathrm{span}(\mathrm{Id} - M_{\overline{G}}) \cap T_{x^\star}^G = \{0\}. \end{aligned}$$

(ii) The proof is classical using the spectral radius formula (2.1), see *e.g.* [3, Theorem 2.12(i)].

(iii) In this case, we have $W_{\overline{G}} = W_{\overline{J}} = \mathrm{Id}$. In turn, $M_{\overline{G}} = \mathrm{P}_{T_{x^\star}^G}$ and $M_{\overline{J}} = \mathrm{P}_{T_{x^\star}^J}$, which yields

$$M = \mathrm{Id} + 2\mathrm{P}_{T_{x^\star}^G}\mathrm{P}_{T_{x^\star}^J} - \mathrm{P}_{T_{x^\star}^G} - \mathrm{P}_{T_{x^\star}^J} = \mathrm{P}_{T_{x^\star}^G}\mathrm{P}_{T_{x^\star}^J} + \mathrm{P}_{S_{x^\star}^G}\mathrm{P}_{S_{x^\star}^J},$$

which is normal, and so is $M_\lambda$. From [2, Proposition 3.6(i)], we get that $\mathrm{Fix}(M) = (T_{x^\star}^J \cap T_{x^\star}^G) \oplus (S_{x^\star}^J \cap S_{x^\star}^G)$. Thus, combining normality, statement (i) and [3, Theorem 2.16] we get that

$$\|M_\lambda^{k+1-K} - M^\infty\| = \|M_\lambda - M^\infty\|^{k+1-K},$$

and $\|M_\lambda - M^\infty\|$ is the optimal convergence rate of $M_\lambda$. Combining together [3, Proposition 3.3] and arguments similar to those of the proof of [2, Theorem 3.10(ii)] (see also [3, Theorem 4.1(ii)]), we get indeed that

$$\|M_\lambda - M^\infty\| = \sqrt{(1-\lambda)^2 + \lambda(2-\lambda)\cos^2\big(\theta_F(T_{x^\star}^J, T_{x^\star}^G)\big)}.$$

The special case is immediate. This concludes the proof. ☐

## C.4 Proof of Corollary 5.5

**Proof.** (i) Let $K \in \mathbb{N}$ sufficiently large such that the locally linearized iteration (5.5) holds. Then we have for $k \geq K$

$$\begin{aligned} z_{k+1} - z^\star &= M_\lambda(z_k - z^\star) + \psi_k + \phi_k = M_\lambda\big(M_\lambda(z_{k-1} - z^\star) + \psi_{k-1} + \phi_{k-1}\big) + \psi_k + \phi_k \\ &= M_\lambda^{k+1-K}(z_K - z^\star) + \textstyle\sum_{j=K}^k M_\lambda^{k-j}(\psi_j + \phi_j). \end{aligned} \tag{C.15}$$

Since $z_k \to z^\star$ from Theorem 3.1 and $M_\lambda$ is convergent to $M^\infty$ by Lemma 5.4(i), taking the limit as $k \to \infty$, we have for all finite $p \geq K$,

$$\lim_{k\to\infty} \sum_{j=p}^{k} M_\lambda^{k-j}(\psi_j + \phi_j) = -M^\infty(z_p - z^\star). \tag{C.16}$$

Using (C.16) in (C.15), we get

$$\begin{aligned}
&z_{k+1} - z^\star \\
&= (M_\lambda - M^\infty)(z_k - z^\star) + \psi_k + \phi_k - \lim_{l\to\infty} \sum_{j=k}^{l} M_\lambda^{l-j}(\psi_j + \phi_j) \\
&= (M_\lambda - M^\infty)(z_k - z^\star) + \psi_k + \phi_k - \lim_{l\to\infty} \sum_{j=k+1}^{l} M_\lambda^{l-j}(\psi_j + \phi_j) - M^\infty(\psi_k + \phi_k) \\
&= (M_\lambda - M^\infty)(z_k - z^\star) + (\mathrm{Id} - M^\infty)(\psi_j + \phi_j) + M^\infty(z_{k+1} - z^\star).
\end{aligned}$$

It is also immediate to see from Lemma 5.4(i) that $\|\mathrm{Id} - M^\infty\| \leq 1$ and

$$(M_\lambda - M^\infty)(\mathrm{Id} - M^\infty) = M_\lambda - M^\infty.$$

Rearranging the terms gives the claimed equivalence.

(ii) Under polyhedrality and constant parameters, we have from Proposition 5.2 that both $\phi_k$ and $\psi_k$ vanish. In this case, (C.16) reads

$$z_k - z^\star \in \ker(M^\infty), \qquad \forall k \geq K,$$

and therefore (5.5) obviously becomes (5.7). □

## C.5 Proof of Theorem 5.6

**Proof.**

(i) Let $K \in \mathbb{N}$ sufficiently large such that (5.6) holds. We then have from Corollary 5.5(i)

$$\begin{aligned}
&(\mathrm{Id} - M^\infty)(z_{k+1} - z^\star) \\
&= (M_\lambda - M^\infty)^{k+1-K}(\mathrm{Id} - M^\infty)(z_K - z^\star) + \sum_{j=K}^{k} (M_\lambda - M^\infty)^{k-j}\big((\mathrm{Id} - M^\infty)\psi_j + \phi_j\big).
\end{aligned}$$

Since $\rho(M_\lambda - M^\infty) < 1$ by Lemma 5.4(i), from the spectral radius formula, we know that for every $\rho \in ]\rho(M_\lambda - M^\infty), 1[$, there is a constant $C$ such that

$$\|(M_\lambda - M^\infty)^j\| \leq C\rho^j$$

for all integers $j$. We thus get

$$\begin{aligned}
&\|(\mathrm{Id} - M^\infty)(z_{k+1} - z^\star)\| \\
&\leq C\rho^{k+1-K}\|z_K - z^\star\| + C\sum_{j=K}^{k} \rho^{k-j}\phi_j + C\sum_{j=K}^{k} \rho^{k-j}\|(\mathrm{Id} - M^\infty)\psi_j\| \\
&= C\rho^{k+1-K}\Big(\|z_K - z^\star\| + \rho^{K-1}\sum_{j=K}^{k} \frac{\phi_j}{\rho^j}\Big) + C\sum_{j=K}^{k} \rho^{k-j}\|(\mathrm{Id} - M^\infty)\psi_j\|,
\end{aligned} \tag{C.17}$$

By assumption, $\phi_j = C'\eta^j$, for some constant $C' \geq 0$ and $\eta < \rho$, and we have

$$\rho^{K-1}\sum_{j=K}^{k} \frac{\phi_j}{\rho^j} \leq C'\rho^{K-1}\sum_{j=K}^{\infty} (\eta/\rho)^j = \frac{C'\eta^K}{\rho - \eta} < +\infty.$$

Setting $C'' = C\left(\|z_K - z^\star\| + \frac{C'\eta^K}{\rho - \eta}\right) < +\infty$, we obtain

$$\|(\mathrm{Id} - M^\infty)(z_{k+1} - z^\star)\| \leq C''\rho^{k+1-K} + C\sum_{j=K}^{k} \rho^{k-j}\|(\mathrm{Id} - M^\infty)\psi_j\|.$$

This, together with the fact that $\|(\mathrm{Id} - M^\infty)\psi_j\| = o(\|(\mathrm{Id} - M^\infty)(z_j - z^\star)\|)$ yields the claimed result.

(ii) From Corollary 5.5(ii), we have

$$z_k - z^\star = (M_\lambda - M^\infty)^{k+1-K}(z_K - z^\star).$$

Moreover, by virtue of Lemma 5.4(iii), $M_\lambda$ is normal and converges linearly to

$$M^\infty = \mathrm{P}_{(T^J_{x^\star} \cap T^G_{x^\star}) \oplus (S^J_{x^\star} \cap S^G_{x^\star})}$$

at the optimal rate

$$\rho = \|M_\lambda - M^\infty\| = \sqrt{(1-\lambda)^2 + \lambda(2-\lambda)\cos^2\left(\theta_F(T^J_{x^\star}, T^G_{x^\star})\right)}.$$

Combining all this then entails

$$\begin{aligned}\|z_{k+1} - z^\star\| &\le \|(M_\lambda - M^\infty)^{k+1-K}\|(z_K - z^\star) \\ &= \|M_\lambda - M^\infty\|^{k+1-K}\|z_K - z^\star\| \\ &= \rho^{k+1-K}\|z_K - z^\star\|.\end{aligned}$$

□

## C.6 Proof of Proposition 7.2

**Proof.**

(i) Similarly to (5.3), we have $W_{\overline{J}}$ is firmly non-expansive. Then for $W_{\overline{G}}$, since $L$ is injective, so is $L_G$, then both $L_G^T L_G$ and $(L_G^T L_G)^{-1}$ are *symmetric positive definite*. Therefore, we have the following similarity result for $W_{\overline{G}}$,

$$\begin{aligned} W_{\overline{G}} &= L_G\Big((L_G^T L_G)^{-\frac{1}{2}}\big(\mathrm{Id} + (L_G^T L_G)^{-\frac{1}{2}} H_{\overline{G}} (L_G^T L_G)^{-\frac{1}{2}}\big)(L_G^T L_G)^{\frac{1}{2}}\Big)^{-1}(L_G^T L_G)^{-1} L_G \\ &= L_G (L_G^T L_G)^{-\frac{1}{2}}\Big(\mathrm{Id} + (L_G^T L_G)^{-\frac{1}{2}} H_{\overline{G}} (L_G^T L_G)^{-\frac{1}{2}}\Big)^{-1}(L_G^T L_G)^{\frac{1}{2}}(L_G^T L_G)^{-1} L_G \\ &= L_G (L_G^T L_G)^{-\frac{1}{2}}\Big(\mathrm{Id} + (L_G^T L_G)^{-\frac{1}{2}} H_{\overline{G}} (L_G^T L_G)^{-\frac{1}{2}}\Big)^{-1}(L_G^T L_G)^{-\frac{1}{2}} L_G. \end{aligned} \tag{C.18}$$

Since $(L_G^T L_G)^{-\frac{1}{2}} H_{\overline{G}} (L_G^T L_G)^{-\frac{1}{2}}$ is symmetric positive definite, then the matrix

$$(\mathrm{Id} + (L_G^T L_G)^{-\frac{1}{2}} H_{\overline{G}} (L_G^T L_G)^{-\frac{1}{2}})^{-1}$$

is firmly non-expansive. It is easy to show that matrix $\|L_G (L_G^T L_G)^{-\frac{1}{2}}\| \le 1$, then owing to [4, Corollary 4.6], $W_{\overline{G}}$ is firmly non-expansive.

(ii) Under the assumptions of Corollary 7.1, the identification theorem of ADMM, there exists $K \in \mathbb{N}$ large enough such that for all $k \ge K$, $(v_k, u_k) \in \mathcal{M}^J_{v^\star} \times \mathcal{M}^G_{u^\star}$. Denote $T^J_{v_k}$ and $T^J_{v^\star}$ be the tangent spaces corresponding to $v_k$ and $v^\star \in \mathcal{M}^J_{v^\star}$, and similarly $T^G_{u_k}$ and $T^G_{u^\star}$ the tangent spaces corresponding to $u_k$ and $u^\star \in \mathcal{M}^G_{u^\star}$. Denote $\tau^J_k : T^J_{v_k} \to T^J_{v^\star}$ (resp. $\tau^G_k : T^G_{u_k} \to T^G_{u^\star}$) the parallel translation along the unique geodesic on $\mathcal{M}^J_{v^\star}$ (resp. $\mathcal{M}^G_{u^\star}$) joining $v_k$ to $v^\star$ (resp. $u_k$ to $u^\star$).
For convenience, define $\gamma' = 1/\gamma$. From (7.2), for $v_k$, we have

$$\begin{cases} v_k = \mathrm{prox}_{\gamma' J}(\gamma' w_k), \\ v^\star = \mathrm{prox}_{\gamma' J}(\gamma' w^\star), \end{cases} \iff \begin{cases} \gamma' w_k - v_k \in \gamma' \partial J(v_k), \\ \gamma' w^\star - v^\star \in \gamma' \partial J(v^\star). \end{cases}$$

Projecting on the corresponding tangent spaces, using Lemma 4.2, and applying the parallel translation operator $\tau^J_k$ leads to

$$\begin{aligned} \gamma' \tau^J_k \nabla_{\mathcal{M}^J_{v^\star}} J(v_k) &= \tau^J_k \mathrm{P}_{T^J_{v_k}}(\gamma' w_k - v_k) \\ &= \mathrm{P}_{T^J_{v^\star}}(\gamma' w_k - v_k) + \big(\tau^J_k \mathrm{P}_{T^J_{v_k}} - \mathrm{P}_{T^J_{v^\star}}\big)(\gamma' w_k - v_k), \\ \gamma' \nabla_{\mathcal{M}^J_{v^\star}} J(v^\star) &= \mathrm{P}_{T^J_{v^\star}}(\gamma' w^\star - v^\star). \end{aligned}$$

We then obtain

$$
\begin{aligned}
&\gamma'\tau_k^J \nabla_{\mathcal{M}^J_{v^\star}} J(v_k) - \gamma' \nabla_{\mathcal{M}^J_{v^\star}} J(v^\star) \\
&= P_{T^J_{v^\star}}\big((\gamma' w_k - \gamma' w^\star) - (v_k - v^\star)\big) \\
&\quad + \underbrace{\big(\tau_k^J \mathrm{P}_{T^J_{v_k}} - \mathrm{P}_{T^J_{v^\star}}\big)(\gamma' w_k - v_k - \gamma' w^\star + v^\star)}_{\textbf{Term 1}} + \underbrace{\big(\tau_k^J \mathrm{P}_{T^J_{v_k}} - \mathrm{P}_{T^J_{v^\star}}\big)(\gamma' w^\star - v^\star)}_{\textbf{Term 2}}.
\end{aligned}
\tag{C.19}
$$

Combining (A.1) and Lemma C.1, we have for **Term 1**

$$
(\tau_k^J \mathrm{P}_{T^J_{v_k}} - \mathrm{P}_{T^J_{v^\star}})(\gamma' w_k - v_k - \gamma' w^\star + v^\star) = o(\gamma' \|w_k - w^\star\|). \tag{C.20}
$$

As far as **Term 2** is concerned, with (7.5), (A.1) and the Riemannian Taylor expansion Lemma C.1, and recall that $\gamma' y^\star = \gamma' w^\star - v^\star$, we have

$$
\begin{aligned}
&\gamma'\tau_k^J \nabla_{\mathcal{M}^J_{v^\star}} J(v_k) - \gamma' \nabla_{\mathcal{M}^J_{v^\star}} J(v^\star) - \big(\tau_k^J \mathrm{P}_{T^J_{v_k}} - \mathrm{P}_{T^J_{v^\star}}\big)(\gamma' w^\star - v^\star) \\
&= \tau_k^J\big(\gamma' \nabla_{\mathcal{M}^J_{v^\star}} J(v_k) - \mathrm{P}_{T^J_{v_k}}(\gamma' w^\star - v^\star)\big) - \big(\gamma' \nabla_{\mathcal{M}^J_{v^\star}} J(v^\star) - \mathrm{P}_{T^J_{v^\star}}(\gamma' w^\star - v^\star)\big) \\
&= \tau_k^J \nabla_{\mathcal{M}^J_{v^\star}} \overline{\overline{J}}(v_k) - \nabla_{\mathcal{M}^J_{v^\star}} \overline{\overline{J}}(v^\star) = P_{T^J_{v^\star}} \nabla^2_{\mathcal{M}^J_{v^\star}} \overline{\overline{J}}(v^\star) P_{T^J_{v^\star}}(v_k - v^\star) + o(\|v_k - v^\star\|) \\
&= P_{T^J_{v^\star}} \nabla^2_{\mathcal{M}^J_{v^\star}} \overline{\overline{J}}(v^\star) P_{T^J_{v^\star}}(v_k - v^\star) + o(\gamma' \|w_k - w^\star\|).
\end{aligned}
\tag{C.21}
$$

Therefore, inserting (C.20) and (C.21) into (C.19), we obtain

$$
\begin{aligned}
&& H_{\overline{\overline{J}}}(v_k - v^\star) &= \gamma' P_{T^J_{v^\star}}(w_k - w^\star) - P_{T^J_{v^\star}}(v_k - v^\star) + o(\gamma' \|w_k - w^\star\|) \\
\Longrightarrow && (\mathrm{Id} + H_{\overline{\overline{J}}}) P_{T^J_{v^\star}}(v_k - v^\star) &= \gamma' P_{T^J_{v^\star}}(w_k - w^\star) + o(\gamma' \|w_k - w^\star\|) \\
\Longrightarrow && P_{T^J_{v^\star}}(v_k - v^\star) &= \gamma' W_{\overline{\overline{J}}} P_{T^J_{v^\star}}(w_k - w^\star) + o(\gamma' \|w_k - w^\star\|) \\
\Longrightarrow && P_{T^J_{v^\star}}(v_k - v^\star) &= \gamma' P_{T^J_{v^\star}} W_{\overline{\overline{J}}} P_{T^J_{v^\star}}(w_k - w^\star) + o(\gamma' \|w_k - w^\star\|) \\
\Longrightarrow && v_k - v^\star &= \gamma' M_{\overline{\overline{J}}}(w_k - w^\star) + o(\gamma' \|w_k - w^\star\|),
\end{aligned}
\tag{C.22}
$$

where we used the fact that $v_k - v^\star = P_{T^J_{v^\star}}(v_k - v^\star) + o(\|v_k - v^\star\|)$ Lemma C.1.
Similarly for $u_{k+1}$, we have

$$
\begin{cases}
L^T(w_k - 2y_k - \gamma L u_{k+1}) \in \partial G(u_{k+1}), \\
L^T(w^\star - 2y^\star - \gamma L u^\star) \in \partial G(u^\star),
\end{cases}
$$

Upon projecting onto the corresponding tangent spaces and applying the parallel translation $\tau^G_{k+1}$, we get

$$
\begin{aligned}
\tau^G_{k+1} \nabla_{\mathcal{M}^G_{u^\star}} G(u_{k+1}) &= \tau^G_{k+1} \mathrm{P}_{T^G_{u_{k+1}}}\big(L^T(w_k - 2y_k - \gamma L u_{k+1})\big) \\
&= \mathrm{P}_{T^G_{u^\star}}\big(L^T(w_k - 2y_k - \gamma L u_{k+1})\big) \\
&\quad + \big(\tau^G_{k+1} \mathrm{P}_{T^G_{u_{k+1}}} - \mathrm{P}_{T^G_{u^\star}}\big)\big(L^T(w_k - 2y_k - \gamma L u_{k+1})\big), \\
\nabla_{\mathcal{M}^G_{u^\star}} G(u^\star) &= \mathrm{P}_{T^G_{u^\star}}\big(L^T(w^\star - 2y^\star - \gamma L u^\star)\big).
\end{aligned}
$$

Subtracting both equations, we obtain

$$\begin{aligned}
&\tau_{k+1}^G \nabla_{\mathcal{M}_{u^\star}^G} G(u_{k+1}) - \nabla_{\mathcal{M}_{u^\star}^G} G(u^\star) \\
&= P_{T_{u^\star}^G}\big(L^T(w_k - 2y_k - \gamma L u_{k+1}) - L^T(w^\star - 2y^\star - \gamma L u^\star)\big) \\
&\quad + \underbrace{\big(\tau_{k+1}^G \mathrm{P}_{T_{u_{k+1}}^G} - \mathrm{P}_{T_{u^\star}^G}\big)\big(L^T(w_k - 2y_k - \gamma L u_{k+1}) - L^T(w^\star - 2y^\star - \gamma L u^\star)\big)}_{\textbf{Term 3}} \\
&\quad + \underbrace{\big(\tau_{k+1}^G \mathrm{P}_{T_{u_{k+1}}^G} - \mathrm{P}_{T_{u^\star}^G}\big)\big(L^T(w^\star - 2y^\star - \gamma L u^\star)\big)}_{\textbf{Term 4}}.
\end{aligned} \tag{C.23}$$

Owing to (7.3), we have

$$\begin{aligned}
\|(y_k + \gamma L u_{k+1}) - (y^\star + \gamma L u^\star)\| &= \|w_{k+1} - w^\star\| \le \|w_k - w^\star\|, \\
\|(w_k - y_k) - (w^\star - y^\star)\| &= \gamma\|v_k - v^\star\| \le \gamma\|\gamma'(w_k - w^\star)\| = \|w_k - w^\star\|.
\end{aligned}$$

Hence, with similar arguments to those used for **Term 1**, we have **Term 3** $= o(\|w_k - w^\star\|)$. Moreover, as $w^\star - 2y^\star - \gamma L u^\star = -y^\star$, then similarly to (C.21), we have for **Term 4**,

$$\begin{aligned}
&\tau_{k+1}^G \nabla_{\mathcal{M}_{u^\star}^G} G(u_{k+1}) - \nabla_{\mathcal{M}_{u^\star}^G} G(u^\star) + \big(\tau_{k+1}^G \mathrm{P}_{T_{u_{k+1}}^G} - \mathrm{P}_{T_{u^\star}^G}\big) L^T y^\star \\
&= \gamma P_{T_{u^\star}^G} \nabla^2_{\mathcal{M}_{u^\star}^G} \overline{\overline{G}}(u^\star) P_{T_{u^\star}^G}(u_{k+1} - u^\star) + o(\|w_k - w^\star\|).
\end{aligned} \tag{C.24}$$

Then for (C.23) we have,

$$\begin{aligned}
H_{\overline{\overline{G}}}(u_{k+1} - u^\star) &= \gamma' P_{T_{u^\star}^G} L^T\big((w_k - 2y_k - \gamma L u_{k+1}) - (w^\star - 2y^\star - \gamma L u^\star)\big) + o(\gamma'\|w_k - w^\star\|) \\
&= \gamma' L_G^T\big((w_k - 2y_k - \gamma L u_{k+1}) - (w^\star - 2y^\star - \gamma L u^\star)\big) + o(\gamma'\|w_k - w^\star\|) \\
&= \gamma' L_G^T(w_k - w^\star) - \gamma' L_G^T(2y_k - 2y^\star) - L_G^T L\big(u_{k+1} - u^\star\big) + o(\gamma'\|w_k - w^\star\|) \\
&= \gamma' L_G^T(w_k - w^\star) - \gamma' L_G^T(2y_k - 2y^\star) - L_G^T L_G(u_{k+1} - u^\star) + o(\gamma'\|w_k - w^\star\|),
\end{aligned}$$

which leads to, (recall that $y_k = w_k - \gamma v_k$)

$$\begin{aligned}
&\big(\mathrm{Id} + (L_G^T L_G)^{-1} H_{\overline{\overline{G}}}\big)(u_{k+1} - u^\star) \\
&= \gamma'(L_G^T L_G)^{-1} L_G^T(w_k - w^\star) - \gamma'(L_G^T L_G)^{-1} L_G^T(2y_k - 2y^\star) + o(\gamma'\|w_k - w^\star\|) \\
&= -\gamma'(L_G^T L_G)^{-1} L_G^T(w_k - w^\star) + 2(L_G^T L_G)^{-1} L_G^T(v_k - v^\star) + o(\gamma'\|w_k - w^\star\|),
\end{aligned}$$

since $V_{\overline{\overline{G}}} = (\mathrm{Id} + (L_G^T L_G)^{-1} H_{\overline{\overline{G}}})^{-1}$, then the above equation can be further reformulated as

$$\begin{aligned}
&\gamma L_G(u_{k+1} - u^\star) \\
&= -L_G V_{\overline{\overline{G}}}(L_G^T L_G)^{-1} L_G^T(w_k - w^\star) + 2\gamma L_G V_{\overline{\overline{G}}}(L_G^T L_G)^{-1} L_G^T(v_k - v^\star) + o(\|w_k - w^\star\|) \\
&= -W_{\overline{\overline{G}}}(w_k - w^\star) + 2\gamma W_{\overline{\overline{G}}}(v_k - v^\star) + o(\|w_k - w^\star\|),
\end{aligned} \tag{C.25}$$

Summing up (C.22) and (C.25), we get

$$\begin{aligned}
w_{k+1} - w^\star &= (w_k - \gamma v_k + \gamma L u_{k+1}) - (w^\star - \gamma v^\star + \gamma L u^\star) \\
&= (w_k - w^\star) - \gamma(v_k - v^\star) + \gamma L(u_{k+1} - u^\star) \\
&= (w_k - w^\star) - M_{\overline{\overline{J}}}(w_k - w^\star) + o(\|w_k - w^\star\|) \\
&\qquad - W_{\overline{\overline{G}}}(w_k - w^\star) + 2\gamma W_{\overline{\overline{G}}}(v_k - v^\star) + o(\|w_k - w^\star\|) \\
&= (w_k - w^\star) - M_{\overline{\overline{J}}}(w_k - w^\star) + o(\|w_k - w^\star\|) - W_{\overline{\overline{G}}}(w_k - w^\star) + 2W_{\overline{\overline{G}}} M_{\overline{\overline{J}}}(w_k - w^\star) \\
&= \big(\mathrm{Id} - M_{\overline{\overline{J}}} - W_{\overline{\overline{G}}} + 2W_{\overline{\overline{G}}} M_{\overline{\overline{J}}}\big)(w_k - w^\star) - M_{\overline{\overline{J}}}(w_k - w^\star) + o(\|w_k - w^\star\|) \\
&= M(w_k - w^\star) + o(\|w_k - w^\star\|).
\end{aligned}$$

(iii) The convergence rate of sequence $\{w_k\}_{k\in\mathbb{N}}$ following the proof of Theorem 5.6. In the following, we simply derive the form of $M$ when both $G$ and $J$ are locally polyhedral around $u^\star$ and $v^\star$ respectively. For this case, $H_{\overline{\overline{G}}}$ and $H_{\overline{\overline{J}}}$ vanish and then $W_{\overline{\overline{G}}}, W_{\overline{\overline{J}}}$ become

$$W_{\overline{\overline{G}}} \stackrel{\text{def}}{=} L_G(L_G^T L_G)^{-1}L_G^T \quad \text{and} \quad W_{\overline{\overline{J}}} \stackrel{\text{def}}{=} P_{T^J_{v^\star}},$$

where $W_{\overline{\overline{G}}}$ now is the projection operator opto the subspace $T^{G,L}_{u^\star}$. As a result, we have

$$M = P_{T^{G,L}_{u^\star}} P_{T^J_{v^\star}} + (\mathrm{Id} - P_{T^{G,L}_{u^\star}})\mathrm{P}_{S^J_{v^\star}}.$$

The optimal convergence result then follows Theorem 5.6.

(iv) When $M^\infty = 0$ for $G, J$ are locally polyhedral around $u^\star, v^\star$ respectively, we have

$$\|w_{k+1} - w^\star\| \le O(\rho^{k-K})$$

hold for any $\rho \in ]\rho(M), 1[$.
We have from (7.2) that $v_k = \mathrm{prox}_{J/\gamma}(w_k/\gamma)$, then owing to the non-expansiveness of $\mathrm{prox}_{J/\gamma}$ (Lemma 2.3), we have

$$\|v_{k+1} - v^\star\| = \|\mathrm{prox}_{J/\gamma}(w_{k+1}/\gamma) - \mathrm{prox}_{J/\gamma}(w^\star/\gamma)\| \le \frac{1}{\gamma}\|w_{k+1} - w^\star\|.$$

Then as $y_{k+1} = \mathrm{prox}_{\gamma J^*}(w_{k+1})$ and $\gamma L u_{k+1} = w_{k+1} - y_k$, we have

$$\begin{aligned} \|y_{k+1} - y^\star\| &= \|\mathrm{prox}_{\gamma J^*}(w_{k+1}) - \mathrm{prox}_{\gamma J^*}(w^\star)\| \le \|w_{k+1} - w^\star\|, \\ \|\gamma L(u_{k+1} - u^\star)\| &= \|(w_{k+1} - y_k) - (w^\star - y^\star)\| \le \|w_{k+1} - w^\star\| + \|y_k - y^\star\|, \end{aligned} \tag{C.26}$$

which leads to the claimed result. □

# References

[1] P.-A. Absil, R. Mahony, and J. Trumpf. An extrinsic look at the Riemannian Hessian. In *Geometric Science of Information*, pages 361–368. Springer, 2013.

[2] H. Bauschke, J.Y.B. Cruz, T.A. Nghia, H.M. Phan, and X. Wang. The rate of linear convergence of the douglas-rachford algorithm for subspaces is the cosine of the friedrichs angle. *J. of Approx. Theo.*, 185(63–79), 2014.

[3] H. H. Bauschke, J.Y. Bello Cruz, T.A. Nghia, H. M. Phan, and X. Wang. Optimal rates of convergence of matrices with applications. *Numerical Algorithms*, 2016. in press (arxiv:1407.0671).

[4] H. H. Bauschke and P. L. Combettes. *Convex analysis and monotone operator theory in Hilbert spaces*. Springer, 2011.

[5] H. H. Bauschke, M. N. Dao, D. Noll, and H. M. Phan. On Slater's condition and finite convergence of the Douglas–Rachford algorithm for solving convex feasibility problems in Euclidean spaces. *Journal of Global Optimization*, pages 1–21, 2015. in press.

[6] H. H. Bauschke and W.M. Moursi. On the order of the operators in the Douglas–Rachford algorithm. *Optimization Letters*, 2016. in press (arXiv:1505.02796v1).

[7] D. Boley. Local linear convergence of the alternating direction method of multipliers on quadratic or linear programs. *SIAM Journal on Optimization*, 23(4):2183–2207, 2013.

[8] J. M. Borwein and B. Sims. The Douglas–Rachford algorithm in the absence of convexity. In H. H. Bauschke, R. S. Burachik, P. L. Combettes, V. Elser, D. R. Luke, and H. Wolkowicz, editors, *Fixed-Point Algorithms for Inverse Problems in Science and Engineering*, volume 49 of *Springer Optimization and Its Applications*, pages 93–109. Springer New York, 2011.

[9] H. Brézis. *Opérateurs maximaux monotones et semi-groupes de contractions dans les espaces de Hilbert*. North-Holland Math. Stud. Elsevier, New York, 1973.

[10] I. Chavel. *Riemannian geometry: a modern introduction*, volume 98. Cambridge University Press, 2006.

[11] P. L. Combettes. Fejér monotonicity in convex optimization. In A. Christodoulos Floudas and M. Panos Pardalos, editors, *Encyclopedia of Optimization*, pages 1016–1024. Springer, Boston, MA, 2001.

[12] P. L. Combettes. Quasi-Fejérian analysis of some optimization algorithms. *Studies in Computational Mathematics*, 8:115–152, 2001.

[13] P. L. Combettes. Solving monotone inclusions via compositions of nonexpansive averaged operators. *Optimization*, 53(5-6):475–504, 2004.

[14] P. L. Combettes and J.-C. Pesquet. A proximal decomposition method for solving convex variational inverse problems. *Inverse Problems*, 24(6):065014, 2008.

[15] P. L. Combettes and I. Yamada. Compositions and convex combinations of averaged nonexpansive operators. *Journal of Mathematical Analysis and Applications*, 425(1):55–70, 2015.

[16] A. Daniilidis, D. Drusvyatskiy, and A. S. Lewis. Orthogonal invariance and identifiability. *SIAM J. Matrix Anal. Appl.*, 35:580–598, 2014.

[17] D. Davis and W. Yin. Convergence rate analysis of several splitting schemes. Technical report, arXiv:1406.4834, August 2014.

[18] D. Davis and W. Yin. Convergence rates of relaxed Peaceman–Rachford and ADMM under regularity assumptions. Technical report, arXiv:1407.5210, 2014.

[19] L. Demanet and X. Zhang. Eventual linear convergence of the Douglas–Rachford iteration for basis pursuit. *Mathematics of Computation*, 85(297):209–238, 2016.

[20] J. Douglas and H. H. Rachford. On the numerical solution of heat conduction problems in two and three space variables. *Transactions of the American mathematical Society*, 82(2):421–439, 1956.

[21] J. Eckstein and D. P. Bertsekas. On the Douglas-Rachford splitting method and the proximal point algorithm for maximal monotone operators. *Mathematical Programming*, 55(1-3):293–318, 1992.

[22] D. Gabay. Applications of the method of multipliers to variational inequalities. In M. Fortin and R. Glowinski, editors, *Augmented Lagrangian Methods: Applications to the Solution of Boundary-Value Problems*, pages 299–331. Elsevier, North-Holland, Amsterdam, 1983.

[23] P. Giselsson and S. Boyd. Metric selection in Douglas–Rachford Splitting and ADMM. *arXiv preprint arXiv:1410.8479*, 2014.

[24] W. L. Hare and A. S. Lewis. Identifying active constraints via partial smoothness and prox-regularity. *Journal of Convex Analysis*, 11(2):251–266, 2004.

[25] W.L. Hare and A. S. Lewis. Identifying active manifolds. *Alg. Op. Res.*, 2(2):75–82, 2007.

[26] R. Hesse and D. R. Luke. Nonconvex notions of regularity and convergence of fundamental algorithms for feasibility problems. *SIAM Journal on Optimization*, 23(4):2397–2419, 2013.

[27] R. Hesse, D. R. Luke, and P. Neumann. Projection methods for sparse affine feasibility: Results and counterexamples. Technical report, 2013.

[28] N.T.B. Kim and D.T. Luc. Normal cones to a polyhedral convex set and generating efficient faces in multiobjective linear programming. *Acta Math. Vietnam.*, 25:101–124, 2000.

[29] J. M. Lee. *Smooth manifolds*. Springer, 2003.

[30] C. Lemaréchal, F. Oustry, and C. Sagastizábal. The U-lagrangian of a convex function. *Trans. Amer. Math. Soc.*, 352(2):711–729, 2000.

[31] A. S. Lewis. Active sets, nonsmoothness, and sensitivity. *SIAM Journal on Optimization*, 13(3):702–725, 2003.

[32] A. S. Lewis, D. R. Luke, and J. Malick. Local linear convergence for alternating and averaged nonconvex projections. *Found. Comput. Math.*, 9(4):485–513, 2009.

[33] J. Liang, M. J. Fadili, and G. Peyré. Local linear convergence of Forward–Backward under partial smoothness. In *Advances in Neural Information Processing Systems*, pages 1970–1978, 2014.

[34] J. Liang, M. J. Fadili, and G. Peyré. Activity identification and local linear convergence of Forward–Backward-type methods, 2015. submitted (arXiv:1503.03703).

[35] J. Liang, M. J. Fadili, and G. Peyré. Convergence rates with inexact non-expansive operators. *Mathematical Programming*, pages 1–32, 2015.

[36] J. Liang, M. J. Fadili, G. Peyré, and R. Luke. Activity identification and local linear convergence of Douglas–Rachford/ADMM under partial smoothness. In *Scale Space and Variational Methods in Computer Vision*, pages 642–653. Springer, 2015.

[37] P. L. Lions and B. Mercier. Splitting algorithms for the sum of two nonlinear operators. *SIAM Journal on Numerical Analysis*, 16(6):964–979, 1979.

[38] F.J. Luque. Asymptotic convergence analysis of the proximal point algorithm. *SIAM Journal on Control and Optimization*, 22:277–293, 1984.

[39] S. A. Miller and J. Malick. Newton methods for nonsmooth convex minimization: connections among-Lagrangian, Riemannian Newton and SQP methods. *Mathematical programming*, 104(2-3):609–633, 2005.

[40] H. M. Phan. Linear convergence of the Douglas-Rachford method for two closed sets. *Optimization*, 65(2):369–385, 2016.

[41] H. Raguet, M. J. Fadili, and G. Peyré. A generalized Forward–Backward splitting. *SIAM Journal on Imaging Sciences*, 6(3):1199–1226, 2013.

[42] R. T. Rockafellar. Monotone operators and the proximal point algorithm. *SIAM Journal on Control and Optimization*, 14(5):877–898, 1976.

[43] R. T. Rockafellar. *Convex analysis*, volume 28. Princeton university press, 1997.

[44] R. T. Rockafellar and R. Wets. *Variational analysis*, volume 317. Springer Verlag, 1998.

[45] S. Vaiter, C. Deledalle, J. M. Fadili, G. Peyré, and C. Dossal. The degrees of freedom of partly smooth regularizers. *Annals of the Institute of Statistical Mathematics*, 2015. to appear.

[46] S. Vaiter, G. Peyré, and M. J. Fadili. Model consistency of partly smooth regularizers. Technical Report arXiv:1307.2342, submitted, 2015.

[47] S. J. Wright. Identifiable surfaces in constrained optimization. *SIAM Journal on Control and Optimization*, 31(4):1063–1079, 1993.